\input amstex
\documentstyle{amsppt}
\document
\topmatter
\title
On the torsion of Drinfeld modules of rank two
\endtitle
\author
Ambrus P\'al
\endauthor
\date
October 22, 2008.
\enddate
\address Department of Mathematics, 180 Queen's Gate, Imperial College,
London SW7 2AZ, United Kingdom\endaddress
\email a.pal\@imperial.ac.uk\endemail
\abstract We prove that the curve $Y_0(\goth p)$ has no $\Bbb F_2(T)$-rational points where $\goth p\triangleleft\Bbb F_2[T]$ is a prime ideal of degree at least $3$ and $Y_0(\goth p)$ is the affine Drinfeld modular curve parameterizing Drinfeld modules of rank two over $\Bbb F_2[T]$ of generic characteristic with Hecke-type level $\goth p$-structure. As a consequence we derive a conjecture of Schweizer describing completely the torsion of Drinfeld modules of rank two over $\Bbb F_2(T)$ implying the uniform boundedness conjecture in this particular case. We reach our results with a variant of the formal immersion method. Moreover we show that the group $\text{\rm Aut}(X_0(\goth p))$ has order two. As a further application of our methods we also determine the prime-to-$p$ cuspidal torsion packet of $X_0(\goth p)$ where $\goth p\triangleleft\Bbb F_q[T]$ is a prime ideal of degree at least $3$ and $q$ is a power of the prime $p$.
\endabstract
\footnote" "{\it 2000 Mathematics Subject Classification. \rm 11G09
(primary), 11G18 (secondary).}
\endtopmatter

\heading 1. Introduction
\endheading

\definition{Notation 1.1} Let $F=\Bbb F_q(T)$ denote the rational function field of transcendence degree one over a finite field $\Bbb F_q$ of characteristic $p$, where $T$ is an indeterminate, and let $A=\Bbb F_q[T]$. Let $Y_0(\goth n)$ denote the Drinfeld modular curve
parameterizing Drinfeld modules of rank two over $\Bbb F_q[T]$ of
generic characteristic with Hecke-type level $\goth n$-structure,
where $\goth n\triangleleft\Bbb F_q[T]$ is a non-zero ideal. Assume now that $q=2$. The main result of this paper is the following
\enddefinition
\proclaim{Theorem 1.2} The curve $Y_0(\goth p)$ has no $F$-rational
points if $\goth p$ is a prime ideal and $\deg(\goth p)\geq3$.
\endproclaim
This result is analogous to the celebrated theorem of Mazur on the
isogenies of elliptic curves of prime degree defined over the field
of rational numbers (proved in [14]). The basic strategy of the latter paper is to prove the formal immersion property at the cusps at every finite place. Mazur's proof of this fact is difficult to adapt to Drinfeld modular curves because of the lack of a good theory of $q$-expansions for Drinfeld modular forms. In fact the latter is unlikely to exist because one of its important consequences, namely multiplicity one, is probably false (see a discussion of this matter in 9.7.4 of [8] on page 90). Instead we will use a variant introduced in the papers [15] and [22] to show the formal immersion property at the special fiber of $Y_0(\goth n)$ over the place $\infty$ corresponding to the point at infinity on $\Bbb P^1_{\Bbb F_q}$. As a consequence in the second half of the proof we are forced to take an approach which is also different from Mazur's original method studying isogeny characters and which is rather particular to function fields (see chapter 9).

It is possible to derive uniform bounds on the torsion of Drinfeld modules of rank two using the theorem above. In order to state this result we need to introduce some notation.
\definition{Notation 1.3} For every $\Bbb F_q$-algebra $B$ let
$B\{\tau\}$ denote the skew-polynomial ring over $B$ defined by the
relation $\tau b=b^q\tau$, where $b$ is any element of $B$. This
ring is naturally isomorphic to the ring of $\Bbb F_q$-linear
endomorphisms of the group scheme $\Bbb G_a$ over $B$ where the
action of $\tau$ is given by the $q$-power map $x\mapsto x^q$ and
the elements of $B$ act by scalar multiplication. Let
$\phi:A\rightarrow B\{\tau\}$ be a Drinfeld module of rank two over
$B$ where $B$ is an $A$-algebra. In this special case $\phi$ is
simply an $\Bbb F_q$-algebra homomorphism with
$\phi(T)=T+g\tau+\Delta\tau^2$ with some $\Delta\in B^*$. Then the
additive group of $B$ becomes an $A$-module via the action of $A$ on
the group scheme $\Bbb G_a$ induced by $\phi$. Let us write
$(^{\phi}B)_{tors}$ for the torsion of this $A$-module. 

When $B$ is actually equal to $F$ it is not difficult to prove that
the group $(^{\phi}F)_{tors}$ is finite, in fact we can always write
it in the form $(^{\phi}F)_{tors}=A/\goth m\oplus A/\goth n$ where
$\goth m$ and $\goth m$ are non-zero ideals of $A$ and $\goth m$
divides $\goth n$. For any ideal $\goth a\triangleleft A$ let
$\deg(\goth a)$ denote the degree of the unique monic polynomial
generating $\goth a$. Assume now that $q=2$. Using the main result of this paper we will prove the following theorem which was conjectured in this precise
form by A$\text{\rm .}$ Schweizer (see Conjecture 1 of [26] on page
601).
\enddefinition
\proclaim{Theorem 1.4} We have $\deg(\goth m)+\deg (\goth
n)\leq2$ for any Drinfeld module $\phi$ of rank two over $F$.
\endproclaim
In fact the above was conjectured by A$\text{\rm .}$ Schweizer for
every $q$ but we are only able to prove it for $q=2$ at the moment.
This theorem is sharp, that is every possibility allowed by the
theorem above actually occurs (see the remark at the bottom of page 604 in [26]). It trivially implies the uniform boundedness conjecture made by B. Poonen in [21] (see page 572) in this particular case.
\definition{Notation 1.5} It is possible to reformulate Theorem 1.4 in an equivalent, elementary form which makes no
reference to Drinfeld modules. Given a polynomial $P(t)\in F[t]$ we
say that $x\in F$ is a preperiodic point for $P(t)$ if the set
$\{P^{[n]}(x)|n\in\Bbb N\}$ is finite, where $P^{[n]}(t)=P\circ
P\circ\cdots\circ P(t)$ is the $n$-fold composition of $P(t)$ with
itself. For a Drinfeld module $\phi:A\rightarrow F\{\tau\}$ the set
of preperiodic points of the additive polynomial $\phi(T)$ in $F$ is
equal to $(^{\phi}F)_{tors}$.
\enddefinition
Assume again that $q=2$.
\proclaim{Theorem 1.6} For every polynomial $P(t)=Tt+gt^2+\Delta t^4\in
F[t]$ the set of preperiodic points of $P(t)$ in $F$ has cardinality
at most $4$.
\endproclaim
\definition{Proof} Note that when $F=\Bbb F_2(T)$ then every Drinfeld module of rank one over $F$ is isomorphic to the Carlitz module. Because the torsion of the latter has cardinality $4$ the theorem above follows immediately from Theorem 1.4.\ $\square$
\enddefinition
We have another application. Let $K$ be an algebraically closed field containing $F$. For every ideal $\goth n$ as above let $X_0(\goth n)$ denote the unique geometrically irreducible smooth projective curve over $F$ which contains $Y_0(\goth n)$ as an open subvariety. When  $\deg(\goth p)\geq 3$ the genus of $X_0(\goth p)$ is at least $2$ hence the group of automorphisms Aut$(X_0(\goth p))$ of the base change of $X_0(\goth p)$ to $K$ is finite and it does not depend on the choice of $K$. Using the main result of this paper (and under the hypothesis $q=2$) we are able to prove the following
\proclaim{Corollary 1.7} The group $\text{\rm Aut}(X_0(\goth p))$ has order two.
\endproclaim
It is interesting to note that the original proof of the analogous theorem (in [16]) uses the $q$-expansion of modular forms in a way which is again difficult to adapt to our situation. We are able to find a different proof by heavily exploiting that $F$ has characteristic $2$. 

In order to prove our main result, Theorem 1.2, we need to use the theory of the Eisenstein ideal developed in [17] and [18] to prove some facts about the Eisenstein quotient of the Jacobian of the Drinfeld modular curve $X_0(\goth p)$, the latter introduced first in [29] by Akio Tamagawa in this context. In particular we answer positively a question of Akio Tamagawa raised at the end of [29]. The latter result can be used to derive the nice application below.
\definition{Notation 1.8} For every ideal $\goth n$ as above let $J_0(\goth n)$ denote the Jacobian of the curve $X_0(\goth n)$. For every $F$-rational point $P_0$ of $X_0(\goth p)$ we may regard $X_0(\goth p)$ as a subvariety of $J_0(\goth p)$ over $F$ via the Albanese embedding $X_0(\goth p)\hookrightarrow J_0(\goth p)$ which maps every closed point $P$ to the class of the divisor $P-P_0$. We let $\Cal S_{P_0}$ denote the prime-to-$p$ torsion packet of $X_0(\goth p)$ relative to $P_0$, that is, the set:
$$\Cal S_{P_0}=\{P\in X_0(\goth p)(\overline F)|\ 
\text{$\exists n\in\Bbb N$ such that $n(P-P_0)=0$ in $J_0(\goth p)(\overline F)$ and $p\!\!\not|n$}\},$$
where for every field $K$ let $\overline K$ denote its maximal separable extension. The geometric points of the zero dimensional complement of $Y_0(\goth n)$ in $X_0(\goth n)$ are called the cusps of the curve $X_0(\goth n)$. When $\goth n=\goth p$ is a prime ideal then the curve $X_0(\goth p)$ has two cusps which are actually defined over $F$. Finally let $\widetilde F$ denote the maximal constant field extension of $F$. 
\enddefinition
\proclaim{Theorem 1.9} Assume that $\deg(\goth p)\geq 3$ and $P_0$ is one of the two cusps of $X_0(\goth p)$. Then we have:
$$\Cal S_{P_0}=\cases X_0(\goth p)(\widetilde F)\cup WP(\goth p),&
\text{when $\deg(\goth p)=3$ and $p\neq2$},\\
X_0(\goth p)(\widetilde F),&
\text{when $\deg(\goth p)=3$ and $p=2$},\\
\{\text{the two cusps}\},&\text{otherwise,}\endcases$$
where $WP(\goth p)$ denotes the set of Weierstrass points of $X_0(\goth p)$. 
\endproclaim
\definition{Contents 1.10} In the next chapter we prove a useful theorem about the cokernel of  maps between the connected components of N\'eron models of certain abelian varieties induced by functoriality, mildly generalizing a result in [19]. This theorem is first put to use in the third chapter. Here we recall the definition of the cuspidal divisor group $\Cal C(\goth p)$ and the Eisenstein quotient of the Jacobian of the Drinfeld modular curve $X_0(\goth p)$, then show that the cuspidal divisor group maps isomorphically onto the Mordell-Weil group of this quotient. Moreover at the end of this chapter we prove Theorem 1.9. It is important to note that as a consequence of the results of this chapter we could use the Eisenstein quotient instead of the winding quotient if we could prove that it also has the crucial formal immersion property. But the period lattice of the Eisenstein quotient, unlike that of the winding quotient, is not explicit, this is why we work with the latter.

In the fourth chapter we recall the definition of modular symbols and some results about these objects which were proved in [18]. In the fifth chapter we construct the winding quotient of the Jacobian of the Drinfeld modular curve $X_0(\goth p)$ when $\goth p$ is a prime ideal of degree at least $3$ by describing explicitly its lattice of periods with respect to its rigid analytic uniformization at $\infty$ and derive its basic properties using this definition. In the sixth chapter we study the structure of quotient of the Bruhat-Tits tree by the Hecke congruence group. Although the results proved in this chapter are completely elementary, they are crucial for the proof. The aim of the following chapter is to construct a model of $X_0(\goth p)$ over the spectrum of $\Cal O_{\infty}$, the valuation ring of the completion $F_{\infty}$ of $F$ with respect to $\infty$. The key properties of this model, namely its regularity and the formal immersion property for the winding quotient, is derived with the help of the results of the previous chapter.

We study the Galois module structure of the mod $\goth p$ torsion of a rank two Drinfeld module over $F$ which has good reduction at every place $v\neq\infty$ of $F$ in the eighth chapter. In the ninth chapter we actually prove Theorem 1.2 and derive Corollary 1.7 as a consequence. First we determine the set of $F_{\infty}$-rational points of $X_0(\goth p)$. When the rational point reduces to one of the cusps we follow the method used in [15] and [22]. In the remaining cases we show that the corresponding Drinfeld module has a constant $j$-invariant and we apply the results of the previous chapter. In the last chapter we prove Theorem 1.4 by analyzing Drinfeld modular curves of small conductor. As a consequence of the work of A$\text{\rm .}$ Schweizer we only need to study the Drinfeld modular curve $X_0(T^4+T^2+1)$ more seriously.
\enddefinition
\definition{Acknowledgment 1.11} I wish to thank Akio Tamagawa for some useful discussions about the proofs of Corollary 1.7 and Theorem 1.9. 
\enddefinition

\heading 2. Component groups of N\'eron models of totally split abelian varieties
\endheading

\definition{Definition 2.1} For every algebraic group $T$ over a field $K$ which is a split torus let $C(T)$ denote its group of cocharacters. Then $C(T)$ is a free and finitely generated abelian group whose rank is equal to the dimension of $T$ over $K$. The group $T(K)$ of $K$-valued points of $T$ is canonically isomorphic to $K^*\otimes C(T)$. Let $K$ be now a field complete with respect to a discrete valuation $v:K^*\rightarrow\Bbb Z$. A subgroup $\Lambda$ of $T(K)=K^*\otimes C(T)$ is called a discrete lattice if the restriction of the homomorphism $v\otimes1:K^*\otimes C(T)\rightarrow
\Bbb Z\otimes C(T)$ to $\Lambda$ is injective and the image has finite cokernel. In this case the quotient $T/\Lambda$ exists in the category of rigid analytic spaces and it is a proper rigid analytic group such that the quotient map $T\rightarrow T/\Lambda$ is a homomorphism of rigid analytic groups. Let $Q$ be another split torus over $K$ and let $\Delta$ be a discrete lattice in $Q(K)=K^*\otimes C(Q)$. Let $\text{Hom}(T,\Gamma;Q,\Delta)$ denote the group of homomorphisms
$\phi:T\rightarrow Q$ of algebraic groups over $K$ such that $\phi(\Gamma)\subseteq\Delta$. The operation of forming quotients induces an injective homomorphism
$$h:\text{Hom}(T,\Gamma;Q,\Delta)\rightarrow\text{Hom}(T/\Gamma,Q/\Delta)$$
where
the latter is the group of rigid analytic homomorphisms from the rigid
analytic group $T/\Lambda$ to the rigid analytic group $Q/\Delta$.
\enddefinition
\proclaim{Theorem 2.2} The homomorphism $h$ is an isomorphism.
\endproclaim
\definition{Proof} This is essentially Satz 5 of [9], page 33, although the
result there is only stated for endomorphisms, i.e$\text{.}$ in the case
when $T$ and $Q$ are equal. One may immediately reduce to this case by
considering $\phi\oplus0:T\oplus Q\rightarrow T\oplus Q$.\ $\square$
\enddefinition
\definition{Notation 2.3} Let $\Cal O$ and $\bold k$ denote the valuation ring of $K$ and its residue field, respectively. For every abelian variety $C$ over $K$ let $\Phi(C)$ denote the group of connected components of the special fiber of the N\'eron model of $C$ over $\Cal O$. Let $A$ be an abelian variety defined over $K$ and assume that the N\'eron model $\Cal A$ of $A$ over $\Cal O$ has split multiplicative reduction. The latter means that the connected component of the identity of the special fiber of $\Cal A$ is a split torus over $\bold k$. Let $B$ be an optimal quotient of $A$. Let $\pi:A\rightarrow B$ be the quotient homomorphism and let $\pi_*:\Phi(A)\rightarrow\Phi(B)$ be the homomorphism induced by $\pi$.
\enddefinition
\proclaim{Proposition 2.4} The cokernel of the homomorphism $\pi_*:\Phi(A)\rightarrow\Phi(B)$ injects into the group $(\bold k^*)^{\dim(A)-\dim(B)}$. In particular its order is annihilated by the order of $\bold k^*$ when $\bold k$ is finite.
\endproclaim
The proposition above mildly generalizes one of the main results of [19], but our proof is much simpler and shorter.
\definition{Proof} By the uniformization theory of Gerritzen the abelian variety $A$ has a rigid analytic uniformization by a split torus $T$ over $K$. Hence the abelian variety $B$, as a quotient of $A$, has a rigid analytic uniformization by a split torus $Q$ over $K$, too. Let $\Gamma\subset T(K)$ and $\Delta\subset Q(K)$ be discrete lattices such that  $A=T/\Gamma$ and $B=Q/\Delta$, respectively, as rigid analytic groups. By Theorem 2.2 there is a unique homomorphism $\phi:T\rightarrow Q$ of algebraic groups over $K$ such that $\phi(\Gamma)\subseteq\Delta$ and the homomorphism induced by $\phi$ via forming quotients is equal to $\pi$. Let $L$ denote the completion of the algebraic closure of $K$ with respect to the unique extension of the valuation of $K$. Because $L$ is an algebraically closed field the homomorphism $A(L)\rightarrow B(L)$ induced by $\pi$ is surjective. Hence $Q(L)=\cup_{d\in\Delta}\phi(T(L))+d$. As a topological space $Q(L)$ is locally homeomorphic to a complete metric space. The set $\Delta$ is countable, hence by Baire's category theorem the closed set $\phi(T(L))$ contains an open ball. Therefore $\phi(T(L))$ is an open subgroup of $B(L)$ so the cokernel of the homomorphism $\phi_*:C(T)\rightarrow C(Q)$ induced by $\phi$ is finite.

Let $S$ denote the scheme-theoretical reduction of the kernel of $\phi$. Then $S$ is a smooth algebraic group and it has finitely many connected components. Let $R$ denote the connected component of the identity of $S$: the algebraic group $R$ is a subtorus of $T$. Its group of cocharacters $C(R)
\subseteq C(T)$ lies in the kernel of the homomorphism $\phi_*$ above. Hence the rank of $C(R)$ is at most $\dim(A)-\dim(B)$. Since the homomorphism $v\otimes1:R(K)=K^*\otimes C(R)\rightarrow
\Bbb Z\otimes C(R)$ restricted to $\Gamma\cap R(K)$ is injective, the rank of $\Gamma\cap R(K)$ is at most $\dim(A)-\dim(B)$, too. The group $\Gamma\cap R(K)$ has finite index in the kernel of the restriction of $\phi:T(K)\rightarrow Q(K)$ to $\Gamma$ hence the rank of the finitely generated $\Bbb Z$-module $\phi(\Gamma)\subset Q(K)$ is at least $\dim(A)$. But $\phi(\Gamma)$ lies in $\Delta$ so the rank of $\phi(\Gamma)$ is exactly $\dim(A)$. In particular it is a lattice in $Q(K)$.

Let $V$ denote the split torus which is the quotient of $T$ by its subtorus $R$. Let $\phi_0:T\rightarrow V$ denote the quotient map and let $\phi_1:V\rightarrow Q$ be the unique homomorphism of algebraic groups such that $\phi=\phi_1\circ\phi_0$. Let $G$ denote the torsion of the group $\phi_0(\Gamma)$. The latter is a finite group of $V(K)$, so it injects into the group $(\bold k^*)^{\dim(A)-\dim(B)}$. Moreover its image with respect to $\phi_1$ is the identity since the latter lies in $\phi(\Gamma)$, which has no torsion. Let $W$ denote the split torus which is the quotient of $V$ by its finite subgroup scheme $G$. Let $\phi_{10}:V\rightarrow W$ denote the quotient map and let $\phi_{11}:W\rightarrow Q$ be the unique homomorphism of algebraic groups such that $\phi_1=\phi_{11}\circ\phi_{10}$. The kernel of the homomorphism $\phi_1$ is a finite group scheme, hence the same holds for $\phi_{11}$ as well. Therefore $\phi_{11}$ maps $\phi_{10}\circ\phi_0(\Gamma)$ into $\phi(\Gamma)$ injectively with finite cokernel. Hence $\phi_{10}\circ\phi_0(\Gamma)$ is a lattice in $W(K)$ by the above.

Let $C$ denote the proper rigid analytic group $W/\phi_{10}\circ\phi_0(\Gamma)$. By construction there are homomorphisms $\pi_0:A\rightarrow C$ and $\pi_1:C\rightarrow B$ of rigid analytic groups such that $\pi=\pi_1\circ\pi_0$. By a Theorem of [10] on page 338 the quotient of a split torus $U$ by a discrete lattice $\Lambda\subset U(K)$ is isomorphic to the rigid analytic variety underlying an abelian variety over $K$ if and only if there is a homomorphism $\sigma$ from $\Lambda$ to the character group Hom$(U,\Bbb G_m)$ of $U$ such that the bilinear map:
$$(\alpha,\beta)\mapsto\sigma(\alpha)(\beta):\Lambda\times\Lambda\rightarrow K^*$$
is symmetric and $v(\sigma(\alpha)(\alpha))>0$ whenever $1\neq\alpha\in\Lambda$. Because $Q/\Delta$ is an abelian variety there is such a homomorphism $\sigma $ from $\Delta$ to Hom$(Q,\Bbb G_m)$. Then the composition of $\phi_{11}:\phi_{10}\circ\phi_0(\Gamma)\rightarrow\Delta$, the homomorphism $\sigma$ and the functorial homomorphism Hom$(Q,\Bbb G_{\goth m})\rightarrow\text{Hom}(W,\Bbb G_m)$ is also a homomorphism of this type from $\phi_{10}\circ\phi_0(\Gamma)$ to $\text{Hom}(W,\Bbb G_m)$. Therefore $C$ is an abelian variety over $K$, hence, by GAGA, the homomorphisms $\pi_0$ and $\pi_1$ are algebraic. Recall that $B$ is an optimal quotient of $A$ if and only if for every factorization of the quotient homomorphism $\pi:A\rightarrow B$ into two homomorphisms $\pi':A\rightarrow C'$ and $\pi'':C'\rightarrow B$ of abelian varieties such that $\dim(C')=\dim(B)$ the homomorphism $\pi''$ is an isomorphism. Hence we get that the map $\pi_1:C
\rightarrow B$ is an isomorphism. By Theorem 2.2 the latter implies that the map $\phi_{11}:W\rightarrow Q$ is also an isomorphism. 

The homomorphism $\phi_0:T(K)\rightarrow V(K)$ is surjective. The cokernel of the homomorphism $\phi_{10}:V(K)\rightarrow W(K)$ is isomorphic non-canonically to the finite group $G$. By the above the map $\phi_{11}:W(K)\rightarrow Q(K)$ is an isomorphism. Hence the 
cokernel of the homomorphism $\phi_{10}:T(K)\rightarrow Q(K)$ is isomorphic to $G$. For every split torus $U$ over $K$ and for every discrete lattice $\Lambda\subset U(K)$ such that the quotient $U/\Lambda$ is an abelian variety over $K$ the composition of the quotient map $U(K)\rightarrow(U/\Lambda)(K)$ and the specialization map $(U/\Lambda)(K)\rightarrow\Phi(U/\Lambda)$ is surjective. Since every quotient of $G$ is a subgroup of $(\bold k^*)^{\dim(A)-\dim(B)}$ the claim is now clear.\ $\square$
\enddefinition
Assume now that $A$ has multiplicative reduction which we do not assume to be split, that is, the connected component of the identity of the special fiber of its N\'eron model is a torus over $\bold k$. Let $B$ be again an optimal quotient of $A$ and let $\pi:A\rightarrow B$, $\pi_*:\Phi(A)\rightarrow\Phi(B)$ be the same as above.
\proclaim{Corollary 2.5} The order of cokernel of the homomorphism $\pi_*:\Phi(A)\rightarrow\Phi(B)$ is relatively prime to the characteristic of $\bold k$.
\endproclaim
\definition{Proof} Over a finite unramified extension $L$ of $K$ the abelian variety $A$ has split multiplicative reduction. Because the formation of N\'eron models commute with smooth base change the group of connected components of the special fiber of the N\'eron models of $A$ and $B$ over the integral closure of $\Cal O$ in $L$ are equal to $\Phi(A)$ and $\Phi(B)$, respectively, and the homomorphism induced by the base change of $\pi$ to $L$ between the group of connected components is equal to $\pi_*$. Hence by Proposition 2.4 the cokernel of the homomorphism $\pi_*:\Phi(A)\rightarrow\Phi(B)$ injects into the group $(\bold l^*)^{\dim(A)-\dim(B)}$ where $\bold l$ denotes the residue field of $L$. Since the multiplicative group of $\bold l$ has no finite subgroup whose order is divisible by the characteristic of $\bold k$, the claim is now obvious.\ $\square$
\enddefinition

\heading 3. The Eisenstein quotient
\endheading

\definition{Notation 3.1} For every $f\in A$ we will let the same symbol denote the ideal generated by $f$ when this will not cause confusion. Similarly for every prime ideal $\goth p\triangleleft A$ we let the same symbol denote the unique valuation of $F$ corresponding to this ideal normalized so that $\goth p(\pi)=1$ for every uniformizing element
$\pi\in F$. There is exactly one non-trivial valuation of $F$ (up to rescaling) that is not of this form which will be denoted by $\infty$ as usual. For every valuation $x$ of $F$ let $F_x$ denote the completion of $F$ with respect to $x$ and let $\Cal O_x$, $\bold f_x$ denote its valuation ring and its residue field, respectively.
\enddefinition
\definition{Definition 3.2} The cyclic group generated by the linear equivalence class of the divisor which is the difference of the two cusps of $X_0(\goth p)$ is called the cuspidal divisor group and it is denoted by $\Cal C(\goth p)$. For every proper ideal $\goth m\triangleleft A$ there is a $\goth m$-th Hecke correspondence on the Drinfeld modular curve $X_0(\goth p)$ which in turn induces an endomorphism of the Jacobian $J_0(\goth p)$ of the curve, called the Hecke operator $T_{\goth m}$ (for a detailed description see for example [8].) Let $\Bbb T(\goth p)$ denote the algebra with unity generated by the endomorphisms $T_{\goth q}$ of the Jacobian $J_0(\goth p)$, where $\goth q\triangleleft A$ is any prime ideal which is relatively prime to $\goth p$. The algebra $\Bbb T(\goth p)$ is known to be commutative. Let $\goth E(\goth p)$ denote the ideal of $\Bbb T(\goth p)$ generated by the elements $T_{\goth q}-q^{\deg(\goth q)}-1$, where $\goth q$ is any prime ideal different from $\goth p$. The algebra $\Bbb T(\goth p)$ will be called the Hecke algebra and $\goth E(\goth p)$ its Eisenstein ideal. They were already introduced in Definition 7.10 of [17] on page 160. Let $\gamma$ denote the ideal $\cap_{n=1}^{\infty}\goth E(\goth p)^n$ and let $E(\goth p)$ denote the quotient of $J_0(\goth p)$ by the smallest abelian sub-variety of $J_0(\goth p)$ which is left invariant under the action of $\Bbb T(\goth p)$ and contains the image of the endomorphism $h:J_0(\goth p)\rightarrow J_0(\goth p)$ for every $h\in\gamma$. The abelian variety $E(\goth p)$ is naturally equipped with an action of $\Bbb T(\goth p)$ which makes the quotient map $J_0(\goth n)\rightarrow E(\goth p)$ a $\Bbb T(\goth p)$-equivariant homomorphism.
\enddefinition
Let ${}^r\Cal E(\goth p)$ and ${}^r\Cal H(\goth p)$ denote the largest torsion subgroup of $J_0(\goth p)(\overline F)$ and $E(\goth p)(\overline F)$ annihilated by $\goth E(\goth p)^r$, respectively.
\proclaim{Lemma 3.3} The quotient homomorphism $\pi:J_0(\goth p)\rightarrow E(\goth p)$ maps ${}^r\Cal E(\goth p)$ isomorphically onto ${}^r\Cal H(\goth p)$.
\endproclaim
The claim above and its proof is just a slight extension of claim $(i)$ of Proposition 4.16 of [29] on page 234.
\definition{Proof} For every prime number $l$, natural number $n$ and Galois module $G$ let $G_l$, $G[l^n]$ denote the maximal $l$-primary submodule of $G$ and the $l^n$-torsion submodule of $G$, respectively. Clearly we only have to show that $\pi$ maps ${}^r\Cal E(\goth p)_l$ isomorphically onto ${}^r\Cal H(\goth p)_l$ for every prime $l$ which divides the order of $\Bbb T(\goth p)/\goth E(\goth p)$. By Theorem 1.2 of [18] we know that $l$ is not equal to $p$. Let $\Bbb T_l(\goth p)$ denote the tensor product $\Bbb T(\goth p)\otimes_{\Bbb Z}\Bbb Z_l$ and let $\gamma_l\triangleleft
\Bbb T_l(\goth p)$ be the ideal generated by $\gamma$. Because $\Bbb T_l(\goth p)$ is finitely generated as a $\Bbb Z_l$-module it has only finitely many maximal ideals $\goth P_1,\goth P_2,\ldots,
\goth P_m$. Let $\Bbb T_{\goth P_k}$ denote the completion of $\Bbb T_l(\goth p)$ with respect to $\goth P_k$ and let $\pi_k:\Bbb T_l(\goth p)\rightarrow\Bbb T_{\goth P_k}$ be the canonical projection for every $k=1,2,\ldots,m$. Then the direct sum $\oplus\pi_k:\Bbb T_l(\goth p)\rightarrow\oplus_k\Bbb T_{\goth P_k}$ of these projections is an isomorphism. Because the composition of the canonical inclusion
$\Bbb Z\rightarrow\Bbb T(\goth p)$ and the quotient map $\Bbb T(\goth p)\rightarrow\Bbb T(\goth p)/\goth E(\goth p)$ is surjective (see Proposition 7.11 of [17] on pages 160-161), there is a unique maximal ideal of $\Bbb T_l(\goth p)$, say $\goth P_1$, containing $\goth E(\goth p)$. Let $\epsilon\in\Bbb T_l(\goth p)$ be the idempotent which corresponds to $(1,0,\ldots,0)\in\oplus_k\Bbb T_{\goth P_k}$. Then the annulator ideal of $\epsilon$ in $\Bbb T_l(\goth p)$ is exactly $\gamma_l$ since the latter is the kernel of the projection $\pi_1$ onto the factor $\Bbb T_{\goth P_1}$ by definition.

For every $n\in\Bbb N$ choose an $\epsilon_n\in\Bbb T(\goth p)$ such that $\epsilon_n\equiv\epsilon\mod l^n$. For every $x\in J_0(\goth p)(\overline F)[l^n]$ and for every $n\leq m\in\Bbb N$ we have $\epsilon_n(x)=\epsilon_m(x)$. Let $\epsilon(x)$ denote this common value: then $x\mapsto\epsilon(x)$ is an endomorphism of $J_0(\goth p)(\overline F)_l$ which commutes with the action of $\Bbb T(\goth p)$. Moreover $\pi(\epsilon(x))=\pi(x)$ because $\pi$ is $\Bbb T(\goth p)$-equivariant, the action of $\Bbb T(\goth p)$ on $E(\goth p)$ factors through $\gamma$ and $\pi_1(\epsilon)=1$ by definition. Let $x\in J_0(\goth p)(\overline F)_l$ be a point such that $\pi(x)=0$. Then there is a finite extension $L$ of $F$, a natural number $k$, elements $a_1,a_2,\ldots,a_k\in\gamma$ and points $y_1,y_2,\ldots,y_k\in J_0(\goth p)(L)$ such that
$$x=a_1(y_1)+a_2(y_2)+\cdots+a_k(y_k).$$
For every sufficiently large $n\in\Bbb N$ we have:
$$\epsilon(x)=\epsilon_n(x)=\epsilon_na_1(y_1)+\epsilon_na_2(y_2)+\cdots+
\epsilon_na_k(y_k).$$
Because $\epsilon_na_j\equiv0\mod l^n$ for every $j=1,2,\ldots,k$ we get that $\epsilon(x)\in l^nJ_0(\goth p)(L)$ for every $n\in\Bbb N$. By the Mordell-Weil theorem $J_0(\goth p)(L)$ is a finitely generated $\Bbb Z$-module hence
$\epsilon(x)$ is torsion and its order is prime to $l$. But is is also $l$-primary torsion hence $\epsilon(x)=0$. We get that $\epsilon(y)$ only depends on $\pi(y)$ for every $y\in J_0(\goth p)(\overline F)_l$. Therefore the map $\epsilon$ induces a $\Bbb T(\goth p)$-equivariant section of the surjective map $\pi:J_0(\goth p)(\overline F)_l \rightarrow E(\goth p)(\overline F)_l$. The existence of such a section guarantees that $\pi$ maps ${}^r\Cal E(\goth p)$ surjectively onto ${}^r\Cal H(\goth p)$. Because $\epsilon\equiv1\mod\goth E(\goth p)^r\Bbb T_l(\goth p)$, we have $\epsilon(x)=x$ for every $x\in{}^r\Cal E(\goth p)_l$ so $\pi$ is injective on ${}^r\Cal E(\goth p)_l$ as well.\ $\square$
\enddefinition
\proclaim{Lemma 3.4} Let $B$ be an optimal quotient of $J_0(\goth p)$ over $F$. Then the Mordell-Weil group $B(F)$ has no $p$-torsion.
\endproclaim
\definition{Proof} By 5.1-5.8 of [5], pages 229-233 the curve $X_0(\goth p)$ has a semistable model over $\Cal O_{\goth p}$ such that the irreducible components of the special fiber are rational curves. Hence it has multiplicative reduction by a classical theorem of Raynaud. Moreover the order of the group of connected components of the N\'eron model of $J_0(\goth p)$ at $\Cal O_{\goth p}$ is relatively prime to $p$ by Lemma 5.9 and Proposition 5.10 of [5], page 234. The abelian variety $B$ is an optimal quotient of $J_0(\goth p)$ hence the order of the group of connected components of its N\'eron model at $\Cal O_{\goth p}$ is also relatively prime to $p$ by Corollary 2.5 above. Therefore it has no $F$-rational $p$-torsion by Lemma 7.13 of [17], page 162.\ $\square$
\enddefinition
\proclaim{Theorem 3.5} The quotient homomorphism $\pi:J_0(\goth p)\rightarrow E(\goth p)$ maps $\Cal C(\goth p)$ isomorphically onto $E(\goth p)(F)$.
\endproclaim
\definition{Proof} By the lemma above $E(\goth p)(F)$ has no $p$-torsion. Also note that the action of $\Bbb T(\goth p)$ on $E(\goth p)$ satisfies the Eichler-Shimura relations because the same is true for the action of $\Bbb T(\goth p)$ on $J_0(\goth p)$ and  $\pi:J_0(\goth p)\rightarrow E(\goth p)$ is equivariant with respect to the action of $\Bbb T(\goth p)$. These two facts can be used to show that the torsion of $E(\goth p)$ is annihilated by $\goth E(\goth p)$ the same way as we proved the same claim for $J_0(\goth p)$ in the argument of Lemma 7.16 of [17], page 163. Hence the torsion subgroup of $E(\goth p)(F)$ lies in $\Cal H(\goth p)$ so it must be equal to the image of $\Cal C(\goth p)$ by Lemma 3.3. The latter maps isomorphically onto its image by the lemma just quoted. Now we only have to show that $E(\goth p)(F)$ is finite.

Let $W_{\goth p}$ denote the Atkin-Lehner involution on the curve $Y_0(\goth p)$ (for its definition and properties see [5].) The latter extends to an involution of $X_0(\goth p)$. The automorphism of $J_0(\goth p)$ induced by this involution will be denoted by $W_{\goth p}$ as well. It commutes with the action of $\Bbb T(\goth p)$ hence it leaves the kernel of the projection $\pi:J_0(\goth p)\rightarrow E(\goth p)$ invariant. By Theorem 5.7 of [29] on page 241 we only have to show that the involution of $E(\goth p)$ induced by $W_{\goth p}$, denoted by the same symbol by the usual abuse of notation, is multiplication by $-1$. Let $l$ be an odd prime dividing the order of $\Bbb T(\goth p)/\goth E(\goth p)$ (such a prime exits by the proof of Proposition 4.14 of [29] on page 233). Because $l$ is automatically different from $p$ it will be sufficient to prove that the action of $W_{\goth p}$ on $E(\goth p)(\overline F)_l$ is multiplication by $-1$. First we are going to prove that
 $E(\goth p)(\overline F)_l=\cup_{r=1}^{\infty}{}^r\Cal H(\goth p)_l$. In fact the kernel of the action of $\Bbb T(\goth p)$ on $E(\goth p)(\overline F)[l^n]$ contains the ideals $\gamma$ and $l^n\Bbb T(\goth p)$. Recall that $\gamma$ is the kernel of the canonical map $\Bbb T(\goth p)\rightarrow\Bbb T_{\goth P}$ where $\Bbb T_{\goth P}$ is the completion of $\Bbb T(\goth p)$ with respect to the unique prime ideal $\goth P\triangleleft\Bbb T(\goth p)$ over $\goth E(\goth p)$. Hence the sum of the ideals $\gamma$ and $l^n\Bbb T(\goth p)$ contains some $r$-th power of $\goth E(\goth p)$, so $\goth E(\goth p)^r$ annihilates $E(\goth p)(\overline F)[l^n]$.

Note that for every involution $w$ acting on a finite abelian group $G$ of odd order the group $G$ decomposes as a direct sum of eigenspaces for $w$ with eigenvalues $-1$ and $1$. Using this fact it is easy to show that if $G$ has a filtration by $w$-invariant subgroups $\{0\}=G_0\subset G_1\subset\cdots\subset G_n=G$ such that the involution on $G_{i+1}/G_i$ induced by $w$ is multiplication by $-1$ for $i=0,1,\ldots,n-1$ then the eigenspace of $w$ for the eigenvalue $1$ is trivial by induction on the length of the filtration. Now we are going to prove by induction on $r$ that the action of $W_{\goth p}$ on ${}^r\Cal H(\goth p)_l$ is multiplication by $-1$. It is sufficient to prove the same for ${}^r\Cal E(\goth p)_l$ by Lemma 3.3. The claim holds for $r=1$ by the proposition below. Now we assume that the claim has been proved for $r$, and let $a_1,a_2,\ldots,a_m$ be a set of elements of $\goth E(\goth p)^r$ such that their class mod $\goth E(\goth p)^{r+1}$ generates the $\Bbb Z$-module $\goth E(\goth p)^r/\goth E(\goth p)^{r+1}$. The map $x\mapsto a_1x\oplus\cdots\oplus a_mx$ is a $W_{\goth p}$-invariant homomorphism ${}^{r+1}\Cal E(\goth p)_l\rightarrow{}^1\Cal E(\goth p)_l^m$ with kernel ${}^r\Cal E(\goth p)_l$, hence the claim holds for $r+1$ as well by the remark above.\ $\square$
\enddefinition
\proclaim{Proposition 3.6} The action of $W_{\goth p}$ on ${}^1\Cal E(\goth p)_l$ is multiplication by $-1$.
\endproclaim
\definition{Proof} Recall that the group ${}^1\Cal E(\goth p)$ was denoted by $\Cal E(\goth p)$ in the paper [17]. By Theorem 2.5 of [20] on page 327 every element of $\Cal E(\goth p)$ is defined over the maximal unramified extension $\widetilde F_{\infty}$ of $F_{\infty}$.
Therefore $\Cal E(\goth p)$ has a filtration:
$$0@>>>\Cal E_0(\goth p)@>>>\Cal E(\goth p)@>>>\Cal E_1(\goth p)@>>>0$$
where $\Cal E_0(\goth p)$ is the subgroup of all elements of $\Cal E(\goth p)$ whose image with respect to the specialization map into the special fiber of the N\'eron model of $J_0(\goth p)$ over the valuation ring of $\widetilde F_{\infty}$ lies in the connected component. The subgroup $\Cal E_0(\goth p)$ is left invariant by $W_{\goth p}$ because the automorphism $W_{\goth p}$ of $J_0(\goth p)$ extends to an involution of the N\'eron model of the Jacobian.

We only need to show that $W_{\goth p}$ acts on
$\Cal E_0(\goth p)_l$ and $\Cal E_1(\goth p)_l$ as multiplication by $-1$. First we are going to show this for the former. Let $l(\goth p)$ denote the largest power of $l$ dividing $N(\goth p)$, the order of the group $\Cal C(\goth p)$. By Theorem 1.2 of [18] the number $N(\goth p)$ is also the index of the Eisenstein ideal $\goth E(\goth p)$ in the Hecke algebra $\Bbb T(\goth p)$. Therefore $N(\goth p)$ annihilates the $l$-primary group $\Cal E_0(\goth p)_l$, so the latter is annihilated by $l(\goth p)$ as well. By Proposition 10.8 on pages 189-190 and Corollary 11.7 on page 194 of [17] we know that the $l$-torsion of $\Cal E_0(\goth p)_l$ has order at most $l$. Hence the group $\Cal E_0(\goth p)_l$ is cyclic and its order divides $l(\goth p)$. On the other hand $\Cal E_0(\goth p)_l$ contains the maximal $l$-primary subgroup of the Shimura subgroup $\Cal S(\goth p)$ (by part $(i)$ of Proposition 8.18 in [17], pages 171-172). As the order of the latter is exactly $l(\goth p)$ (by Lemma 8.17 of [17] on page 171) the groups
$\Cal E_0(\goth p)_l$ and $\Cal S(\goth p)_l$ are equal. Because this argument works for every prime $l$ dividing $N(\goth p)$ we get that actually the groups $\Cal E_0(\goth p)$ and $\Cal S(\goth p)$ are equal.

By part $(ii)$ of Proposition 8.18 in [17] on pages 171-172 the group $\Cal S(\goth p)$ maps isomorphically onto the group of connected components of the fiber of the N\'eron model of $J_0(\goth p)$ at $\Cal O_{\goth p}$ via the specialization map. The same is true for the cuspidal divisor group $\Cal C(\goth p)$ by 5.11 of [5] on page 235. Both maps are equivariant with respect to the action of the Atkin-Lehner operator. The latter interchanges the two cusps of $X_0(\goth p)$ hence it act as multiplication by $-1$ on the group $\Cal C(\goth p)$, therefore on the groups $\Cal S(\goth p)$ and $\Cal E(\goth p)_l$ as well.

The group $\Cal E_1(\goth p)_l$ injects into the $l$-primary group of connected components of the special fiber of the N\'eron model of $J_0(\goth p)$ over the valuation ring of  $\widetilde F_{\infty}$ annihilated by the Eisenstein ideal. The latter is a cyclic group by the strong multiplicity one theorem (see the proof of Proposition 7.18 of [17] on pages 163-164) so the same holds for $\Cal E_1(\goth p)_l$ as well. Every involution acting on a cyclic group of odd prime power order is multiplication by $\pm1$ hence it will be sufficient to prove that the action of $W_{\goth p}$ on the $l$-torsion $\Cal E_1(\goth p)[l]$ of $\Cal E_1(\goth p)_l$ is multiplication by $-1$. If $l$ does not divide $q-1$ then $\Cal C(\goth p)_l$ maps isomorphically onto the group $\Cal E_1(\goth p)[l]$ by Proposition 10.8 of [17] on pages 189-190, so the claim above is clear. If $l$ does divide $q-1$ then the group $\Cal D(\goth p)[l]$ of order $l^2$ introduced in Definition 9.17 of [17] on page 183 maps surjectively onto the group $\Cal E_1(\goth p)[l]$ and the kernel is the $l$-torsion of the Shimura group by claim $(v)$ Proposition 9.18 of [17] on pages 184-185. The short
  exact sequence of $l$-torsion Galois modules:
$$0@>>>\Cal S(\goth p)[l]@>>>\Cal D(\goth p)[l]@>>>{\Cal D(\goth p)\over\Cal S(\goth p)[l]}@>>>0$$
is not split by claim $(v)$ of the proposition just quoted above. Hence the involution $W_{\goth p}$ must be multiplication by $-1$ on $\Cal D(\goth p)[l]$ because $W_{\goth p}$ commutes with the Galois action and it is multiplication by $-1$ on $\Cal S(\goth p)[l]$.\ $\square$
\enddefinition
\proclaim{Corollary 3.7} The action of $W_{\goth p}$ on $\Cal E(\goth p)$ is multiplication by $-1$. 
\endproclaim
\definition{Proof} By Lemma 3.3 the \'etale group scheme $\Cal E(\goth p)$ imbeds $W_{\goth p}$-equivariantly into the Eisenstein quotient $E(\goth p)$. On the other hand by the proof of Theorem 3.5 above we know that $W_{\goth p}$ acts on $E(\goth p)$ as multiplication by $-1$.\ $\square$
\enddefinition
Recall that in this chapter we do not assume that $q=2$. 
\proclaim{Proposition 3.8} The set $X_0(\goth p)(\widetilde F)$ lies in
$\Cal S_{P_0}$ when $\goth p$ is a prime ideal of degree $3$ and $P_0$ is a cusp.
\endproclaim
\definition{Proof} Let $X_{\overline F}$ denote the base change of every algebraic variety $X$ defined over $F$ to $\overline F$. Fix a prime $l$ different from $p$ and consider the $l$-adic representation
$H^1(J_0(\goth n)_{\overline F},\overline{\Bbb Q}_l)$ of the
absolute Galois group of $F$ where $\goth n$ is any non-zero ideal
of $A$. By Drinfeld's fundamental theorem (Theorem 2 of [2] on page 562) the latter decomposes as a
sum of irreducible two-dimensional $l$-adic Galois representations
over $\overline{\Bbb Q}_l$. The conductor of every such
representation $\rho$ is the sum of the divisors $\goth m$ and
$\infty$ where $\goth m$ divides $\goth n$. Moreover $\rho$ appears
as an irreducible component of the Galois representation
$H^1(J_0(\goth m)_{\overline F}, \overline{\Bbb Q}_l)$ when $\goth
m$ is a proper divisor of $\goth n$. The curve $X_0(1)$ has genus zero hence the degree of the conductor of every irreducible representation $\rho$ appearing in
$H^1(J_0(\goth p)_{\overline F},\overline{\Bbb Q}_l)$ is $4$.
Therefore the degree of the Grothendieck $L$-function $L(\rho_n,t)$ of the base change $\rho_n$ of $\rho$ to the field $F_n=\Bbb F_{q^n}
(T)$ as a polynomial in $t$ is zero by the Grothendieck-Ogg-Shafarevich formula for every $n\in\Bbb N$. In particular $L(\rho_n,t)$ does not vanish at $q^{-n}$. We get that the Mordell-Weil group $J_0(\goth p)(F_n)$ is finite by the main theorem of [24] on page 509. Hence the group $J_0(\goth p)(\widetilde F)$ is torsion. Therefore it is equal to $\Cal E(\goth p)$ by Theorem 2.5 of [20] on page 327. The order of $\Cal E(\goth p)$ is prime to the characteristic, hence the claim is now clear.\ $\square$
\enddefinition
\definition{Proof of Theorem 1.9} Let $P_0$ and $P_{\infty}$ denote the two cusps of $X_0(\goth p)$. Because the group $\Cal C(\goth p)$ is finite and its order is prime to $p$ both cusps are elements of $\Cal S_{P_0}$. Moreover the Drinfeld modular curve $X_0(\goth p)$ is hyperelliptic when $\deg(\goth p)=3$ by Theorem 20 of [25] on page 343.  By Theorem 20 of [25] quoted above the unique hyperelliptic involution of $X_0(\goth p)$ is the Atkin-Lehner involution $W_{\goth p}$. Hence every Weierstrass point of $X_0(\goth p)$ is a branch point of the unique hyperelliptic covering (see Proposition 1.4 of [30] on page 286). The Atkin-Lehner involution exchanges $P_0$ and $P_{\infty}$ so for every hyperelliptic branch point $P$ of $X_0(\goth p)$ we have:
$$2P\sim P_{\infty}+P_0\quad\text{or equivalently}\quad
2(P-P_0)\sim P_{\infty}-P_0,$$
where $\sim$ denotes the linear equivalence of divisors. In particular $P$ is a $2N(\goth p)$-torsion point with respect to the Albanese embedding $X_0(\goth p)\hookrightarrow J_0(\goth p)$ with base point $P_0$. Hence it is an element of $\Cal S_{P_0}$ when $p$ is different from $2$.

Now let $P$ be an element of $\Cal S_{P_0}$ and assume that it has order $n\in\Bbb N$. By assumption $n$ is relatively prime to $p$. Let $G$ denote the Galois module generated by $P$ in $J_0(\goth p)(\overline F)[n]$. Then $G$ is unramified at every place of $F$ different from $\goth p$ and $\infty$ because at these places the abelian variety $J_0(\goth p)$ has good reduction. At the places $\goth p$ and $\infty$ the curve $J_0(\goth p)$ has stable reduction. Hence the Galois module $G$ is either also unramified at the places $\goth p$ and $\infty$ or the curve $X_0(\goth p)$ is hyperelliptic and $P$ is a Weierstrass point by Proposition 0.2 of [30] on page 283. The latter is only possible when $\deg(\goth p)=3$ by Theorem 20 of [25] quoted above and the order of $P-P_0$ is $2N(\goth p)$ in this case as we already saw above. Hence we may assume that the Galois module $G$ is everywhere unramified. In this case $P$ is defined over the maximal constant field extension $\widetilde F$ of $F$. By Theorem 2.5 of [20] on page 327 the torsion of the group $J_0(\goth p)(\widetilde F)$ is equal to $\Cal E(\goth p)$. Hence $P-P_0$ lies in
$\Cal E(\goth p)$, so by Corollary 3.7 we have:
$$W_{\goth p}(P-P_0)\sim P_0-P\quad\text{or equivalently}\quad
P+W_{\goth p}(P)\sim P_0+P_{\infty}.$$
If $P$ is different from $P_0$ and $P_{\infty}$ the relation above implies that $X_0(\goth p)$ is hyperelliptic. In the latter case $\deg(\goth p)=3$, hence the claim follows from Proposition 3.8. \ $\square$
\enddefinition

\heading 4. The winding homomorphism
\endheading

\definition{Definition 4.1} For any graph $G$ let $\Cal V(G)$ and $\Cal
E(G)$ denote its set of vertices and edges, respectively. In this paper
we will only consider such oriented graphs $G$ which are equipped with an involution $\overline{\cdot}:\Cal E(G)\rightarrow\Cal E(G)$ such that for each edge $e\in\Cal E(G)$ the original and terminal vertices of the edge $\overline e\in\Cal E (G)$ are the terminal and original vertices of $e$, respectively. The edge $\overline e$ is called the edge $e$ with reversed orientation. Let $R$ be a commutative group. A function $\phi:\Cal E(G)\rightarrow R$ is called a harmonic $R$-valued cochain, if it satisfies the following conditions:
\roster
\item"$(i)$" We have:
$$\phi(e)+\phi(\overline e)=0\text{\ }(\forall e\in\Cal E(G)).$$
\item"$(ii)$" If for an edge $e$ we introduce the notation $o(e)$ and
$t(e)$ for its original and terminal vertex respectively,
$$\sum_{{e\in\Cal E(G)\atop o(e)=v}}\phi(e)=0\text{\ }
(\forall v\in\Cal V(G)).$$
\endroster
We denote by $H(G,R)$ the group of $R$-valued harmonic cochains on $G$.
\enddefinition
\definition{Definition 4.2} Let $GL_2$ denote the group scheme of invertible two by two matrices and let $Z$ denote its center. Let $\pi\in F_{\infty}$ be a uniformizer. We are going to recall the definition of the Bruhat-Tits tree $\Cal T$ associated to the projective linear group $PGL_2(F_{\infty})$. The set of vertices $\Cal V(\Cal T)$ and edges $\Cal E(\Cal T)$ are the cosets $GL_2(F_{\infty})/GL_2(\Cal O_{\infty})Z(F_{\infty})$ and $GL_2(F_{\infty})/\Gamma_{\infty}Z(F_{\infty})$, respectively, where $\Gamma_{\infty}$ is the Iwahori group:
$$\Gamma_{\infty}=\left\{\left(\matrix a&b\\
c&d\endmatrix\right)\in GL_2(\Cal O_{\infty})|
\infty(c)>0\right\}.$$
Since $\Gamma_{\infty}$ is a subgroup of $GL_2(\Cal O_{\infty})$ there is a natural map $o:\Cal E(\Cal T)\rightarrow\Cal V(\Cal T)$ which assigns to every edge its original vertex. The matrix $\left(
\smallmatrix0&1\\\pi&0\endsmallmatrix\right)$ normalizes the Iwahori subgroup therefore the map $GL_2(F_{\infty})\rightarrow GL_2(F_{\infty})$
given by the rule $\left(\smallmatrix a&b\\c&d\endsmallmatrix\right)
\rightarrow\left(\smallmatrix a&b\\c&d\endsmallmatrix\right)
\left(\smallmatrix0&1\\\pi&0\endsmallmatrix\right)$
induces a map on the coset $\Cal V(\Cal T)$. This map is the involution which assigns to every edge $e$ the same edge $\overline e$ with reversed orientation. The composition of this involution and the map $o$ is the map $t:
\Cal E(\Cal T)\rightarrow\Cal V(\Cal T)$ which assigns to every edge its terminal vertex.
\enddefinition
\definition{Definition 4.3} For every non-zero ideal $\goth n
\triangleleft A$ let $\Gamma_0(\goth n)$ denote the Hecke congruence group:
$$\Gamma_0(\goth n)=\{\left(\matrix a&b\\c&d\endmatrix\right)
\in GL_2(A)|c\equiv0\text{\ mod }\goth n\}.$$ The group
$GL_2(F_{\infty})$ acts on itself via its left-regular action which
induces an action of $GL_2(F_{\infty})$ on the Bruhat-Tits tree.
This action induces an action of its subgroup $\Gamma_0(\goth n)$ on
$\Cal T$ as well. Let $H(\Cal T,R)^{\Gamma_0(\goth n)}$ denote the
group of $\Gamma_0(\goth n)$-invariant $R$-valued cochains on $\Cal
T$. The group $GL_2(A)$ does not contain elements which map an edge
$e\in\Cal E(\Cal T)$ to the same edge $\overline e$ with reversed
orientation therefore the cosets $\Gamma_0(\goth n)\backslash\Cal
V(\Cal T)$ and $\Gamma_0(\goth n)\backslash\Cal E(\Cal T)$ are the
vertices and edges of an oriented graph which is going to be denoted
by $\Gamma_0(\goth n)\backslash\Cal T$. Every element $\phi$ of
$H(\Cal T,R)^{\Gamma_0(\goth n)}$ induces an $R$-valued function on
the edges of $\Gamma_0(\goth n)\backslash\Cal T$. If this function
is zero outside of a finite set we say that $\phi$ is cuspidal. The
$R$-module of cuspidal elements of $H(\Cal T,R)^{\Gamma_0(\goth n)}$
is denoted by $H_!(\Cal T,R)^{\Gamma_0(\goth n)}$.
\enddefinition
\definition{Definition 4.4} For every pair $\goth m$, $\goth n
\triangleleft A$ of non-zero ideals let $H(\goth m,\goth n)$ denote the set:
$$H(\goth m,\goth n)=
\{\left(\matrix a&b\\c&d\endmatrix\right)\in GL_2(F)|a,b,c,d\in A,(ad-cb)=\goth m,\goth n\supseteq(c),(d)+\goth n=A\}.$$
The set $H(\goth m,\goth n)$ is a double $\Gamma_0
(\goth n)$-coset and it is also the disjoint union of finitely many left
$\Gamma_0(\goth n)$-cosets. Let $R(\goth m,\goth n)$ be a set of
representatives of these cosets. For any left $\Gamma_0(\goth n)$-invariant $R$-valued function $\phi:\Cal E(\Cal T)\rightarrow R$ define $T_{\goth m}(\phi)$ by the formula:
$$T_{\goth m}(\phi)(g)=\sum_{h\in R(\goth m,\goth n)}\phi(hg),\quad\forall g\in\Cal E(\Cal T).$$
It is well-known and easy to check that $T_{\goth m}(\phi)$ is independent of the choice of $R(\goth m,\goth n)$ and it is also a left $\Gamma_0(\goth n)$-invariant $R$-valued function so we have an $R$-linear operator $T_{\goth m}$ acting on the $R$-module of left $\Gamma_0(\goth n)$-invariant $R$-valued functions on $\Cal E(\Cal T)$. It is also well-known and not too difficult to verify that $T_{\goth m}$ leaves the submodules $H(\Cal T,R)^{\Gamma_0(\goth n)}$ and $H_!(\Cal T,R)^{\Gamma_0(\goth n)}$ invariant. The operator
$T_{\goth m}$ is denoted by the same symbol we use for the operators introduced in Definition 3.2, but this will not cause confusion as we will see. For the moment it is sufficient to remark that they act
on different objects.
\enddefinition
\definition{Definition 4.5} A path $\gamma$ on an oriented
graph $G$ is a sequence of edges $$\{\ldots,e_1,e_2,\ldots,e_n,\ldots\}\in\Cal E(G)$$
indexed by the set $I$ where $I=\Bbb Z$, $I=\Bbb N$ or $I=\{0,1,\ldots,m\}$ for some $m\in\Bbb N$ such that $t(e_i)=o(e_{i+1})$ for every $i$, $i+1\in I$. We say that $\gamma$ is an infinite path, a half-infinite path or a finite path whether we are in the first, in the second or in the third case, respectively. For each edge $e\in\Cal E(G)$ let $i_e:\Cal E(G)\rightarrow\Bbb Z$ denote the unique function such that
$$i_e(f)=
\cases+1,&\text{if $f=e$,}\\
-1,&\text{if $f=\overline e$,}\\
\quad\!\!0,&\text{otherwise.}\endcases$$ Let $\gamma$ be a path
$\{\ldots,e_1,\ldots,e_n,\ldots\}$ on $G$ such that every edge in
$\Cal E(G)$ is only listed finitely many times in the sequence
above. Then the function $i_{\gamma}=\sum_{j\in\Bbb Z}i_{e_j}$ is
well-defined as the sum above has only finitely many terms non-zero
on $e$ for every edge $e\in\Cal E(G)$. Let us consider now the
special case $G=\Gamma\backslash\Cal T$ where $\Gamma=\Gamma_0(\goth
n)$ is a short-hand notation introduced for convenience. Let
$z(\Gamma)$ denote the cardinality of the center of $\Gamma$ and
$\Gamma_e$ is the stabilizer of the edge $e\in\Cal E(\Cal T)$ in
$\Gamma$. (It is well-known that the latter is finite.) For every
path $\gamma$ on the graph $\Gamma\backslash\Cal T$ such that
$i_{\gamma}$ is defined in the sense above we define the function
$\gamma^*:\Cal E(\Cal T)\rightarrow\Bbb Z$ given by the rule
$\gamma^*(e)=|\Gamma_e|i_{\gamma}(\widetilde e)/z(\Gamma)$, where
$\widetilde e$ is the image of the edge $e$ in $\Cal E(\Gamma
\backslash\Cal T)$ and the absolute sign $|\cdot|$ denotes the
cardinality of every finite set. (Since the center of $\Gamma$
leaves the Bruhat-Tits tree invariant, it lies in the stabilizer
$\Gamma_e$, therefore the expression above is indeed an integer.)
\enddefinition
\definition{Definition 4.6} Next we are going to define the fundamental arch connecting two different points $a$, $b\in\Bbb P^1(F_{\infty})$ on the Bruhat-Tits tree. We say that a path $\{\ldots,e_1,\ldots,e_n,\ldots\}$ indexed by the set $I$ on an oriented graph $G$ is without backtracking if $\overline{e_i}\neq e_{i+1}$ for every $i$, $i+1\in I$. Let $S(a,b)$ denote the set of those edges of $\Cal T$ which can be represented by a matrix
$\left(\smallmatrix\alpha&\beta\\\gamma&\delta\endsmallmatrix\right)$ such that the homogeneous coordinates $(\alpha:\gamma)=a$ and $(\beta:\delta)=b$. The elements of the set $S(a,b)$ can be indexed uniquely by the set of integers such that it becomes an infinite path without backtracking: this is the fundamental arch $\overline{ab}$ connecting $a$ and $b$. Let $\overline{ab}$ denote the image of the fundamental arch under the canonical map $\Cal T\rightarrow\Gamma_0(\goth n)\backslash\Cal T$ as well by slight abuse of notation. Let $[a,b]:\Cal E(\Cal T)\rightarrow\Bbb Z$ denote the function $(\overline{ab})^*$ introduced in Definition 4.5 if the latter is well-defined.
\enddefinition
Let $\goth p\triangleleft A$ be now a proper non-zero prime ideal.
\proclaim{Proposition 4.7} The following holds:
\roster
\item"$(i)$" for every different $a$, $b\in\Bbb P^1(F)$ the function $[a,b]$ is well-defined and it is a $\Bbb Z$-valued left ${\Gamma_0(\goth p)}$-invariant harmonic cochain,
\item"$(ii)$" we have $[a,b]\in H_!(\Cal T,\Bbb Z)^{\Gamma_0(\goth p)}$ for every different $a$, $b\in\Bbb P^1(F)$ which are equivalent under the M\"obius action of $\Gamma_0(\goth p)$,
\item"$(iii)$" for every proper non-zero prime ideal $\goth p\neq\goth q
\triangleleft  A$ we have
$$(1+q^{\deg(\goth q)}-T_{\goth q})[0,\infty]\in
(q-1)H_!(\Cal T,\Bbb Z)^{\Gamma_0(\goth p)}.$$
\endroster
\endproclaim
\definition{Proof} This is Proposition 5.3 of [18].\ $\square$
\enddefinition
Let  $r$ be the unique monic polynomial generating the prime ideal $\goth q$. Let $R(\goth q)\subset\Bbb F_q[T]$ denote the set of non-zero polynomials whose degree is less than $\deg(\goth q)$. The following lemma was shown during the proof of Proposition 2.7 in [18].
\proclaim{Lemma 4.8} We have:
$$(1+q^{\deg(\goth q)}-T_{\goth q})[0,\infty]=\sum_{0\neq a\in R(\goth q)}[0,a/r]
\text{.\ $\square$}$$
\endproclaim
\definition{Definition 4.9} Let $\widehat{\Bbb T}(\goth p)$,
$\Bbb T(\goth p)$ denote the commutative $\Bbb Z$-algebra with unity
generated by the endomorphisms $T_{\goth q}$ of the $\Bbb Z$-module $H(\Cal T,\Bbb Z)^{\Gamma_0(\goth p)}$ and $H_!(\Cal T,\Bbb
Z)^{\Gamma_0(\goth p)}$, respectively, where $\goth q \triangleleft
A$ is any prime ideal different from $\goth p$. Note that the latter algebra is denoted by the same symbol as the Hecke algebra introduced in Definition 3.2. Indeed the two algebras are naturally isomorphic by the Gekeler-Reversat uniformization theory (see for example Theorem 7.9 of [17] on page 159 for a convenient description of this isomorphism). By Corollary 3.13 of
[18] the algebras $\widehat{\Bbb T}(\goth p)$ and $\Bbb T(\goth p)$
are sub-algebras of the endomorphism ring of a finitely generated,
free $\Bbb Z$-module hence they must be finitely generated, free
$\Bbb Z$-modules, too. Clearly $\Bbb T(\goth p)$ is the quotient of
$\widehat{\Bbb T}(\goth p)$. Let $\widehat{\goth E}(\goth p)$ denote
the ideal of $\widehat{\Bbb T}(\goth p)$ generated by the elements
$T_{\goth q}-q^{\deg(\goth q)}-1$, where $\goth q\neq\goth p$ is any
prime. We have a well-defined homomorphism
$$e:\widehat{\goth E}(\goth p)
\rightarrow H_!(\Cal T,\Bbb Z)^{\Gamma_0(\goth p)}$$ of
$\widehat{\Bbb T}(\goth p)$-modules given by the rule
$\alpha\mapsto\alpha([0,\infty])/(q-1)$ according to part $(iii)$ of
Proposition 4.7. This map is the analogue of the winding
homomorphism introduced by Mazur.
\enddefinition
The following claim is an immediate corollary of Theorem 5.13 of
[18].
\proclaim{Proposition 4.10} The image of the winding
homomorphism $e$ is non-zero.\ $\square$
\endproclaim

\heading 5. The winding quotient
\endheading

\definition{Definition 5.1} Recall that a finite path
$\{e_0,e_2,\ldots,e_n\}\in\Cal E(G)$ on an oriented graph $G$ is
closed if the equality $t(e_n)=o(e_0)$ holds, too. We define
$H_1(G,\Bbb Z)$ as the abelian group of $\Bbb Z$-valued functions on
$\Cal E(G)$ generated by the functions $i_{\gamma}$ where $\gamma$
is a closed path. We define the map
$$j_{\Gamma_0(\goth n)}:H_1(\Gamma_0(\goth n)\backslash\Cal T,\Bbb Z)
\rightarrow H_!(\Cal T,\Bbb Z)^{\Gamma_0(\goth n)},$$ as the unique
homomorphism which maps $i_{\gamma}$ to the cochain $\gamma^*$ for every $\gamma$ closed path, using the notations of Definition 4.5.
It is easy to see that the homomorphism is well-defined, that is
$\gamma^*$ is indeed a harmonic cochain. By a theorem of Gekeler and Nonnengardt (Theorem 3.3 of [7], page 702) this homomorphism is in fact an isomorphism.
\enddefinition
\definition{Definition 5.2} Let
$\Gamma_0(\goth p)_{\text{ab}}=\Gamma_0(\goth p)/[\Gamma_0(\goth
p),\Gamma_0(\goth p)]$ be the abelianization of $\Gamma_0(\goth p)$,
and let $\overline{\Gamma}_0(\goth p)=\Gamma_0(\goth p)_{\text{ab}}/
(\Gamma_0(\goth p)_{\text{ab}})_{\text{tors}}$ be its maximal
torsion-free quotient. For each $\gamma\in\Gamma_0(\goth p)$ let
$\overline{\gamma}$ denote its image in $\overline{\Gamma}_0(\goth
p)$. Fix a vertex $v_0\in\Cal V(\Cal T)$ and for every
$\gamma\in\Gamma_0(\goth p)$ let $e_0,e_1,\ldots,e_{n(\gamma)}$ be
the unique geodesic path connecting $v_0$ with $\gamma(v_0)$, that
is $v_0=o(e_1)$ and $\gamma(v_0)=t(e_{n(\gamma)})$. Recall that a
path is geodesic if it is the shortest connecting its endpoints, in
this case $v_0$ with $\gamma(v_0)$, i.e$\text{.}$ the number
$n(\gamma)$ is the smallest possible. The image of the path
$e_0,e_1,\ldots$ is closed in $\Gamma_0(\goth p) \backslash\Cal T$:
let $i(\gamma)$ denote the corresponding element in
$H_1(\Gamma_0(\goth p)\backslash\Cal T,\Bbb Z)$. The function $i$
induces a homomorphism $i:\overline{\Gamma}_0(\goth p)\rightarrow
H_1(\Gamma_0(\goth p)\backslash\Cal T,\Bbb Z)$ which is independent
of the choice of $v_0$ and it is an isomorphism. We will use this
identification without further notice. In particular we equip $\overline{\Gamma}_0(\goth p)$ with an action of the Hecke operator $T_{\goth
q}$ via the isomorphism $i$.
\enddefinition
\definition{Notation 5.3} Let $\Bbb C_{\infty}$ denote the completion
of the algebraic closure of $F_{\infty}$ and let
$c: \overline{\Gamma}_0(\goth p)\rightarrow\text{\rm Hom}
(\overline{\Gamma}_0(\goth p),\Bbb C_{\infty}^*)$ be the period map
defined in Proposition 7.5 of [17] on page 159. Note that for every 
finitely generated, free $\Bbb Z$-module $\Cal Y$ the group $\text{\rm Hom}(\Cal Y,\Bbb C_{\infty}^*)$ is canonically isomorphic to the group of $\Bbb C_{\infty}$-valued points of the unique torus $T_{\Cal Y}$ whose group of characters Hom$(T_{\Cal Y},\Bbb G_m)$ is $\Cal Y$. The image of $\overline{\Gamma}_0(\goth p)$ under $c$ is a discrete lattice in 
$\text{\rm Hom}(\overline{\Gamma}_0(\goth p),\Bbb C_{\infty}^*)$, and the rigid-analytic group variety $ T_{\overline{\Gamma}_0(\goth p)}/c(\overline{\Gamma}_0(\goth p))$ is isomorphic to the rigid-analytic variety associated to $J_0(\goth p)$ via the Abel-Jacobi map of Gekeler-Reversat. (For a description of this map using the
same notation see sections 7.1-7.7 of [17] on pages 158-159.)
\enddefinition
\definition{Definition 5.4} Let $\infty:\Bbb
C^*_{\infty}\rightarrow\Bbb Q$ denote a valuation of $\Bbb
C_{\infty}$ which induces its standard non-archimedean topology. Let
$\Cal Y$ be a finitely generated, free $\Bbb Z$-module and let $\infty:
\text{\rm Hom}(\Cal Y,\Bbb C_{\infty}^*)\rightarrow \text{\rm
Hom}(\Cal Y,\Bbb Q)$ denote also the homomorphism induced by
$\infty$. We will say that a submodule $\Lambda$ of $\text{\rm Hom}(\Cal
Y,\Bbb C_{\infty}^*)$ is a quasi-lattice if
\roster
\item"$(i)$" it is a finitely generated abelian group,
\item"$(ii)$" the kernel of the restriction of $\infty$ to $\Lambda$
is finite,
\item"$(iii)$" the free abelian group $\infty(\Lambda)$ spans $\text{\rm
Hom}(\Cal Y,\Bbb Q)$  as a vector space over $\Bbb Q$,
\item"$(iv)$" the $\Bbb Z$-rank of $\infty(\Lambda)$ is the same as the
$\Bbb Z$-rank of $\Cal Y$.
\endroster
When $\Lambda$ is a quasi-lattice in $\text{\rm Hom}(\Cal Y,\Bbb C_{\infty}^*)$ then the quotient $T_{\Cal Y}/\Lambda$ exits in the category of rigid analytic spaces, and it is a proper rigid analytic group such that the quotient map $T_{\Cal Y}\rightarrow T_{\Cal Y}/\Lambda$ is a homomorphism of rigid analytic groups. 
\enddefinition
Our next aim is to give a convenient description of optimal
quotients of $J_0(\goth p)$, considered as an abelian variety over
$\Bbb C_{\infty}$, in terms of certain period lattices. We will say, following [19], that a subgroup $\Lambda$ of $\overline{\Gamma}_0(\goth p)$ is saturated if for every $x\in \overline{\Gamma}_0(\goth p)$ we have $x\in\Lambda$ when
$nx\in\Lambda$ for some positive integer $n$. 
\proclaim{Proposition 5.5} Let $\Cal Y$ be a saturated $\Bbb T(\goth
p)$-invariant subgroup of $\overline{\Gamma}_0(\goth p)$. Then the image $\Lambda$ of $\Cal Y$ under the composition:
$$\CD\overline{\Gamma}_0(\goth p)@>c>>\text{\rm Hom}
(\overline{\Gamma}_0(\goth p),\Bbb C_{\infty}^*)
@>>>\text{\rm Hom}
(\Cal Y,\Bbb C_{\infty}^*)\endCD$$
is a quasi-lattice in $\text{\rm Hom}(\Cal Y,\Bbb C_{\infty}^*)$, and the quotient $ T_{\Cal Y}/\Lambda$ is the rigid-analytic variety attached to a unique optimal quotient of $J_0(\goth p)$. Moreover every optimal quotient of $J_0(\goth p)$ arises this way. 
\endproclaim
\definition{Proof} This is the content of Propositions 3.11 and  3.12 of [19] on pages 2183-2184.\ $\square$
\enddefinition
For every abelian variety $A$ let End$(A)$ denote the ring of endomorphisms of $A$ over its field of definition. The following lemma will be useful.
\proclaim{Lemma 5.6} The following holds:
\roster
\item"$(i)$" Every endomorphism of $J_0(\goth p)$ is defined
over $F$.
\item"$(ii)$" The algebra $\text{\rm End}(J_0(\goth p))\otimes\Bbb Q$ is the product of totally real number fields.
\item"$(iii)$" Every finite subgroup of the group of automorphisms of the algebraic group $J_0(\goth p)$ has exponent $2$.
\endroster
\endproclaim
\definition{Proof} Claim $(i)$ is just Theorem A.1 of [29] on page 242 whose proof is on pages 244-245. Claim $(ii)$ is the content of Remark A.6 of [29] on page 243. The group of automorphisms of the algebraic group $J_0(\goth p)$ is a subgroup of the multiplicative group of the algebra $\text{\rm End}(J_0(\goth p))\otimes\Bbb Q$. Every finite subgroup of the latter has exponent $2$ by the second claim because the maximal finite multiplicative subgroup of a totally real field is $\{\pm1\}$. The last claim is now clear.\ $\square$
\enddefinition
\proclaim{Proposition 5.7} Every optimal quotient of $J_0(\goth p)$ is defined over $F$.
\endproclaim
\definition{Proof} For every abelian variety $A$ let $A^{\vee}$ denote its dual. Let $B$ an optimal quotient of $J_0(\goth p)$ and let $\pi:J_0(\goth p)\rightarrow B$ be the quotient map. Let $A$ be the group scheme which is the kernel of $\pi$. By our assumptions it is an abelian variety. It will be sufficient to prove that $A$ is defined over $F$. Let $\pi^{\vee}:B^{\vee}\rightarrow J_0(\goth p)^{\vee}$ be the dual of $\pi$. Because $B$ is an optimal quotient of $J_0(\goth p)$ the map $\pi^{\vee}$ is an immersion. Since $J_0(\goth p)$ is a Jacobian it has a principal polarization $\iota:J_0(\goth p)\rightarrow
J_0(\goth p)^{\vee}$. Let $\lambda:B\rightarrow B^{\vee}$ be a polarization of $B$. Then the composition
$h=\iota^{-1}\circ\pi^{\vee}\circ\lambda\circ\pi$ is an endomorphism of $J_0(\goth
p)$ hence it is defined over $F$ by claim $(i)$ of the lemma above. In particular
its kernel is a group scheme defined over $F$. Because $\lambda$ is an isogeny the reduction of the connected component of the kernel of $h$ is $A$. Therefore
the latter is also defined over $F$.\ $\square$
\enddefinition
\definition{Notation 5.8} Let $\phi\in H_!(\Cal T,\Bbb C)^{\Gamma_0(\goth p)}$ and $\psi\in H(\Cal T,\Bbb C)^{\Gamma_0(\goth p)}$ be two harmonic cochains. In this case their Petersson product is well-defined (for definition see 4.8 of [8] on page 57) and it will be denoted by
$\langle\phi,\psi\rangle$. Moreover for every $\phi\in H_!(\Cal T,\Bbb C)^{\Gamma_0(\goth p)}$ let $L(\phi,s)$ denote its $L$-function. (For the definition of the latter see [33].)
\enddefinition
We are going to use the following well-known fact (see Proposition 1 of [31] on page 112):
\proclaim{Proposition 5.9} For every $\phi\in H_!(\Cal T,\Bbb C)^{\Gamma_0(\goth p)}$ we have:
$$\langle\phi,[0,\infty]\rangle=c\cdot L(\phi,q^{-1})$$
where $c\in\Bbb C^*$ only depends on $q$.\ $\square$
\endproclaim
\definition{Definition 5.10} Let $\goth I(\goth p)$ denote the ideal:
$$\goth I(\goth p)=\{T\in\widehat{\Bbb T}(\goth p)|T([0,\infty])\in
 H_!(\Cal T,\Bbb Z)^{\Gamma_0(\goth p)}\}$$
of the ring $\widehat{\Bbb T}(\goth p)$. Let $\Lambda(\goth p)$ denote the smallest saturated subgroup of $\overline{\Gamma}_0(\goth p)$ containing the group $\{T([0,\infty])|T\in\goth I(\goth p)\}$. Because the latter is $\Bbb T(\goth p)$-invariant so is the former. Note that
$\Lambda(\goth p)$ contains the image of the winding homomorphism $e$ introduced in Definition 4.9. Let $W(\goth p)$ denote the optimal quotient of $J_0(\goth p)$ corresponding to $\Lambda(\goth p)$ under the construction described in Proposition 5.5.
\enddefinition
For every abelian variety $A$ defined over $F$ let $L(A,t)$ denote its Hasse-Weil
$L$-function. 
\proclaim{Proposition 5.11} The
following holds:
\roster
\item"$(i)$" the quotient $W(\goth p)$ is non-trivial,
\item"$(ii)$" the quotient $W(\goth p)$ is actually defined over $F$,
\item"$(iii)$"  the quotient $W(\goth p)$ is the largest optimal quotient
of $J_0(\goth p)$ such that the Hasse-Weil $L$-function of $W(\goth p)$ does
not vanish at $q^{-1}$.
\endroster
\endproclaim
\definition{Proof} By Proposition 4.10 the group $\Lambda(\goth p)$ is non-zero so claim $(i)$ is clear. The second claim follows from Proposition 5.7. The abelian variety $W(\goth p)$ is isogenous to the product of its absolutely simple optimal quotients, which are all defined over $F$ by Proposition 5.7, since they are also optimal quotients of $J_0(\goth p)$. Hence it will be sufficient to prove that for every absolutely irreducible optimal quotient $B$
of $J_0(\goth p)$ we have $L(B,q^{-1})\neq0$ if and only if $B$ is also an
optimal quotient of $W(\goth p)$. Now let $B$ be an absolutely simple optimal quotient of $J_0(\goth p)$ and let $\Cal Y\subseteq\overline{\Gamma}_0(\goth p)$ be the corresponding saturated $\Bbb T(\goth p)$-invariant subgroup. Note that $B$ is an optimal quotient of $W(\goth p)$ if and only if $\Cal Y\subseteq\Lambda(\goth p)$. Because $B$ is absolutely irreducible
the $\Bbb T(\goth p)\otimes\Bbb Q$-module
$$\Cal Y\otimes\Bbb Q\subseteq\Cal H_0(\Cal T,\Bbb Z)^{\Gamma_0(\goth p)}\otimes
\Bbb Q=\Cal H_0(\Cal T,\Bbb Q)^{\Gamma_0(\goth p)}$$
must be irreducible. Therefore the image of $\Bbb T(\goth p)\otimes\Bbb Q$ in the endomorphism ring of the vector space $\Cal Y\otimes\Bbb Q$ is a number field $K$ and $\Cal Y\otimes\Bbb Q$ has dimension one with respect to this $K$-action. In particular there is a non-zero Hecke eigenform $\phi\in\Cal Y\otimes K\subseteq\Cal H_0(\Cal T,K)^{\Gamma_0(\goth p)}$ such that we have $\Cal Y\subseteq\Lambda(\goth p)$ if and only if $\phi\in\Lambda(\goth p)\otimes K$. Choose an embedding of $K$ into $\Bbb C$. After this choice we may consider $\phi$ as an element of $\Cal H_0(\Cal T,\Bbb C)^{\Gamma_0(\goth p)}$. The natural action of $\Bbb T(\goth p)$ on the optimal quotient $B$ furnishes a homomorphism $K\rightarrow\text{End}(B)$. By the Eichler-Shimura isomorphism (see 8.3.8 of [8] on page 81) the $L$-series $Z_B(q^{-s};K)$ of the abelian variety $B$ with coefficients in $K$ (in the sense of [29] on page 214) is equal to $c_{\phi}L(\phi,s)$ where $c_{\phi}\in\Bbb C$ is a non-zero scalar. Because $Z_B(q^{-s};K)$ vanishes at $s=1$ exactly when $L(B,q^{-1})$ is zero, we only have to show that for every non-zero Hecke eigenform $\phi\in\Cal H_0(\Cal T,\Bbb C)^{\Gamma_0(\goth p)}$ we have $L(\phi,q^{-1})\neq0$ if and only if $\phi\in\Lambda(\goth p)
\otimes\Bbb C$. First assume the latter. Because $\phi$ is an eigenform there is a character $\chi:\widehat{\Bbb T}(\goth p)\otimes\Bbb C\rightarrow\Bbb C$ such that $T(\phi)=\chi(T)\phi$ for every $T\in\widehat{\Bbb T}(\goth p)\otimes\Bbb C$. By assumption there is an operator $T\in\widehat{\Bbb T}(\goth p)\otimes\Bbb C$ such that $\phi=T([0,\infty])$. Because the Petersson product is positive definite we get that:
$$0<\langle\phi,\phi\rangle=\langle\phi,T([0,\infty]\rangle=
\langle T(\phi),[0,\infty]\rangle=\chi(T)\langle\phi,[0,\infty]\rangle=
c\chi(T)L(\phi,q^{-1}),$$
where we also used that the action of $\widehat{\Bbb T}(\goth p)\otimes\Bbb C$ is self-adjoint with respect to the Petersson product in the second equation. Hence $L(\phi,q^{-1})$ is non-zero as we claimed. Assume now that $\phi$ is not an element of $\Lambda(\goth p)
\otimes\Bbb C$. Because the algebra $\widehat{\Bbb T}(\goth p)\otimes\Bbb C$ is semi-simple the vector space $\Lambda(\goth p)
\otimes\Bbb C$ is spanned by Hecke eigenforms. None of these is a scalar multiple of $\phi$ hence $\phi$ is orthogonal to every element of $\Lambda(\goth p)\otimes\Bbb C$. In particular for every non-zero prime ideal $\goth q$ of $A$ different from $\goth p$ we have:
$$\split0=\langle\phi,(T_{\goth q}-(q^{\deg(\goth q)}+1))([0,\infty])\rangle=&\langle(T_{\goth q}-(q^{\deg(\goth q)}+1))(\phi),[0,\infty]\rangle\\=&
(\chi(T_{\goth q})-(q^{\deg(\goth q)}+1))\langle\phi,[0,\infty]\rangle\\
=&c(\chi(T_{\goth q})-(q^{\deg(\goth q)}+1))L(\phi,q^{-1}),\endsplit$$
where we used claim $(iii)$ of Proposition 4.7 in the first equation and the self-adjointness of $T_{\goth q}$ in the second equation. By the Ramanujan-Petersson conjecture (proved in [2] first in this case) the number $\chi(T_{\goth q})$ is not equal to
$1+q^{\deg(\goth q)}$ when $\deg(\goth q)$ is sufficiently large
hence $L(\phi,q^{-1})$ must be zero.\ $\square$
\enddefinition
\proclaim{Corollary 5.12} The following holds:
\roster
\item"$(i)$" the Mordell-Weil group $W(\goth p)(F)$ is finite,
\item"$(ii)$" the homomorphism $J_0(\goth p)(F)\rightarrow
W(\goth p)(F)$ induced by the quotient map is injective restricted to
$\Cal C(\goth p)$.
\endroster
\endproclaim
\definition{Proof} The first claim follows at once from claim $(iii)$ of
Proposition 5.11 above by the main theorem of [24] on page 509. As we saw in the proof of Theorem 3.5 we have $L(E(\goth p),q^{-1})\neq0$ by Theorem 5.7 of [29] on page 241. Hence the Eisenstein quotient $E(\goth p)$ is an optimal quotient of $W(\goth p)$ by claim $(iii)$ of Proposition 5.11. Hence the quotient homomorphism $J_0(\goth p)\rightarrow E(\goth p)$ factors through the quotient homomorphism $J_0(\goth p)\rightarrow W(\goth p)$. Therefore
the induced homomorphism $J_0(\goth p)(F)\rightarrow W(\goth p)(F)$ must be injective restricted to $\Cal C(\goth p)$ by Theorem 3.5.\ $\square$
\enddefinition

\heading 6. The structure of the graph $\Gamma_0(\goth p)
\backslash\Cal T$
\endheading

Some results of this chapter can be found in the paper [7] but they are not stated sufficiently explicitly for our purposes. The computations
of the above quoted paper are also quite involved because they deal with a much more general situation. Hence for the sake of exposition we will use essentially the same methods to derive the necessary results.
\definition{Notation 6.1} Let $v_n$ denote the vertex of the Bruhat-Tits tree represented by the matrix $\left(\smallmatrix T^n&0\\0&1
\endsmallmatrix\right)$ for every natural number $n\in\Bbb N$. Let $G_n$ denote the stabilizer of the vertex $v_n$ in $GL_2(A)$ for every $n$ as above. By Proposition 3 of 1.6 in [28] on pages 86-87 we have:
$$G_0=GL_2(\Bbb F_q),$$
and
$$G_n=
\left\{ \left(\matrix a&b\\0&c\endmatrix\right)\in GL_2(A)|\deg(b)\leq n\right\}$$
when $n>0$. Moreover the vertices $v_n$ form a fun\-da\-men\-tal domain for the action of $GL_2(A)$ on the set of ver\-tices of the
Bruhat-Tits tree and the quotient graph $GL_2(A)\backslash\Cal T$ is an infinite tree such that the image of $v_0$ has degree one and every other vertex has degree two, by the corollary following the proposition just quoted above. Here the degree of a vertex $w$ is the number of edges with origin $w$. 
\enddefinition
\definition{Definition 6.2} Let $u_n$ denote the vertex
$\left(\smallmatrix0&1\\ 1&0\endsmallmatrix\right)
\cdot v_{n+\deg(\goth p)}$ of the Bruhat-Tits tree for every natural number $n\in\Bbb N\cup\{-1\}$. Assume  now that $\deg(\goth p)$ is even and $q=2$. Let $t$ be one of the two elements of $A$ with the following properties:
\roster
\item"$(i)$" we have $\deg(t)<\deg(\goth p)$,
\item"$(ii)$" the reduction of $t$ modulo $\goth p$ generates a degree two extension over $\Bbb F_2\subset\bold f_{\goth p}$. 
\endroster
Let $w$ denote the vertex
$\left(\smallmatrix0&1\\1&t\endsmallmatrix\right)\cdot v_0$
of the Bruhat-Tits tree. For the sake of simple notation let
$\Gamma_v$ denote the stabilizer of $v$ in $\Gamma_0(\goth p)$ for every vertex $v$ of $\Cal T$. 
\enddefinition
\proclaim{Proposition 6.3} Assume that $q=2$ an let $v$ be a vertex of $\Cal T$. Then the stabilizer $\Gamma_v$ of $v$ is non-trivial if and only if $v$ is equivalent to one of the vertices $u_n$, $v_n$ (for some $n\in\Bbb N$) or $w$ modulo the action of $\Gamma_0(\goth p)$. 
\endproclaim
In order to prove the claim above we will need two lemmas.
\proclaim{Lemma 6.4} The set
$$M(\goth p)=\left\{
\left(\matrix0&1\\1&u\endmatrix\right)|u\in A,\deg(u)<\deg(\goth p)
\right\}
\cup\left\{\left(\matrix1&0\\0&1\endmatrix\right)\right\}$$
is a complete set of representatives of left $\Gamma_0(\goth p)$-cosets in $GL_2(A)$.
\endproclaim
\definition{Proof} Let $B$ denote the group scheme of upper-triangular two-by-two matrices and let $GL_2^0(\bold f_{\goth p})$ and $B^0(\bold f_{\goth p})$ denote the subgroup of those elements in $GL_2(\bold f_{\goth p})$ and $B(\bold f_{\goth p})$, respectively, whose determinant lies in $\Bbb F_q^*\subseteq\bold f_{\goth p}^*$. Then we have:
$$[GL_2(A):\Gamma_0(\goth p)]=[GL_2^0(\bold f_{\goth p}):
B^0(\bold f_{\goth p})]=[GL_2(\bold f_{\goth p}):B(\bold f_{\goth p})]=
q^{\deg(\goth p)}+1,$$
because the determinant map restricted to $B(\bold f_{\goth p})$ is a surjection onto $\bold f_{\goth p}^*$. Hence it will be sufficient to prove that the reductions of the elements of $M(\goth p)$ modulo $\goth p$ represent different left $B(\bold f_{\goth p})$-cosets. It is clear that for every $u\in A$ the matrices $\left(\smallmatrix0&1\\1&u\endsmallmatrix\right)$ and $\left(\smallmatrix1&0\\0&1\endsmallmatrix\right)$ are not in the same coset since the latter lies in $B(\bold f_{\goth p})$, but the former does not. Now let $u$, $v\in A$ and assume that there are $a$, $b\in A^*$ and $c\in A$ such that
$$\left(\matrix0&1\\1&u\endmatrix\right)=
\left(\matrix a&b\\0&c\endmatrix\right)\cdot
\left(\matrix0&1\\1&v\endmatrix\right)=
\left(\matrix b&a+bv\\c&cv\endmatrix\right)\mod\goth p.$$
Then $c\equiv1\mod\goth p$ hence $u\equiv v\mod\goth p$. Therefore we must have $u=v$ because $\deg(u)$ and $\deg(v)$ are both less than $\deg(\goth p)$.\ $\square$
\enddefinition
Now we are going to assume that $q=2$.
\proclaim{Lemma 6.5} Let $n$ be a positive integer and let $a$, $b\in A$ such that $\deg(a)<\deg(\goth p)$ and $\deg(b)<\deg(\goth p)$. Then the vertices $\left(\smallmatrix 0&1\\ 1&a\endsmallmatrix\right)
\cdot v_n$ and $\left(\smallmatrix 0&1\\ 1&b\endsmallmatrix\right)\cdot
 v_n$ are equivalent modulo the action of $\Gamma_0(\goth p)$ if and only if $\deg(a-b)\leq n$.
\endproclaim
\definition{Proof} Assume first that $\deg(a-b)\leq n$. Then
$$\left(\matrix0&1\\1&a\endmatrix\right)^{-1}\cdot
\left(\matrix0&1\\1&b\endmatrix\right)=
\left(\matrix -a&1\\1&0\endmatrix\right)\cdot
\left(\matrix0&1\\1&b\endmatrix\right)=
\left(\matrix1&b-a\\ 0&1
\endmatrix\right)\in G_n$$
hence the vertices $\left(\smallmatrix 0&1\\ 1&a\endsmallmatrix\right)
\cdot v_n$ and $\left(\smallmatrix 0&1\\ 1&b\endsmallmatrix\right)
\cdot v_n$ are in fact equal in this case. Assume now that the vertices $\left(\smallmatrix 0&1\\ 1&a\endsmallmatrix\right)
\cdot v_n$ and $\left(\smallmatrix 0&1\\ 1&b\endsmallmatrix\right)\cdot
 v_n$ are equivalent modulo the action of $\Gamma_0(\goth p)$. Then there is an $h\in\Gamma_0(\goth p)$ and $\left(\smallmatrix1&u\\ 0&1\endsmallmatrix\right)\in G_n$ such that 
$$\left(\matrix0&1\\1&a\endmatrix\right)^{-1}\!\!\!\!\cdot h\cdot
\left(\matrix0&1\\ 1&b\endmatrix\right)=
\left(\matrix1&u\\0&1\endmatrix\right).$$
By definition the reduction of $h$ modulo $\goth p$ is an upper-triangular matrix, hence there are elements $k$, $l$, and
$m\in\bold f_{\goth p}$ such that
$$\split\left(\matrix0&1\\1&a\endmatrix\right)^{-1}\!\!\!\!\cdot
\left(\matrix k&l\\ 0&m\endmatrix\right)\cdot
\left(\matrix -a&1\\1&0\endmatrix\right)\equiv&
\left(\matrix -al+m&-ak-alb+mb\\ l&k+lb\endmatrix\right)\\ \equiv&
\left(\matrix1& u\\ 0&1\endmatrix\right)\mod\goth p.\endsplit$$
Therefore we must have $l=0$, $k=m=1$, and $u\equiv b-a\mod\goth p$. Because $\deg(u)\leq n$ the claim follows.\ $\square$
\enddefinition
\definition{Proof of Proposition 6.3} By Lemma 6.4 the vertex $v$ can be written in the form $h\cdot g\cdot v_n$ for some $g\in M(\goth p)$, $h\in\Gamma_0(\goth p)$ and $n\in\Bbb N$. Because we have
$\Gamma_v=h(\Gamma_{g\cdot v_n})h^{-1}$ we may assume that $h=1$. When $g$ is the identity matrix then $v$ is equal to the vertex $v_n$. In this case the stabilizer $\Gamma_v$ is equal to the intersection $\Gamma_0(\goth p)\cap G_n$. The latter is non-trivial for every $n\in\Bbb N$. Hence we may assume that
$g=\left(\smallmatrix0&1\\1&u\endsmallmatrix\right)$ for some $u\in A$ such that $\deg(u)<\deg(\goth p)$. Assume first that $n>0$. For every $b\in A$ we have:
$$\left(\matrix0&1\\1&u\endmatrix\right)\cdot
\left(\matrix1&b\\0&1\endmatrix\right)\cdot
\left(\matrix0&1\\1&u\endmatrix\right)^{-1}=
\left(\matrix0&1\\1&b+u\endmatrix\right)\cdot
\left(\matrix-u&1\\1&0\endmatrix\right)=
\left(\matrix1&0\\ b&1\endmatrix\right).$$
In particular the matrix above is an element of $\Gamma_0(\goth p)$ if and only if $b\in\goth p$. Hence the stabilizer 
$\Gamma_v=g(G_n)g^{-1}\cap\Gamma_0(\goth p)$ is non-trivial if and only if $n\geq\deg(\goth p)$. Moreover by Lemma 6.5 the vertex $v=\left(\smallmatrix0&1\\1&u\endsmallmatrix\right)\cdot v_n$ is equivalent to $u_{n-\deg(\goth p)}$ modulo the action of $\Gamma_0(\goth p)$ in this case. Assume now that $n=0$. For every
$\left(\smallmatrix a&b\\ c&d\endsmallmatrix\right)\in GL_2(\Bbb F_2)$ we have:
$$\split\left(\matrix0&1\\1&u\endmatrix\right)\cdot
\left(\matrix a&b\\ c&d\endmatrix\right)\cdot
\left(\matrix0&1\\1&u\endmatrix\right)^{-1}=&
\left(\matrix c&d\\ a+uc&b+ud\endmatrix\right)\cdot
\left(\matrix-u&1\\1&0\endmatrix\right)\\=&
\left(\matrix -uc+d&c\\ -u^2c+u(d-a)+b&a+uc\endmatrix\right).
\endsplit$$
The lower left entry $-u^2c+u(d-a)+b$ most vanish modulo $\goth p$ if the product matrix above is an element of $\Gamma_0(\goth p)$.  Hence the reduction of $u$ modulo $\goth p$ generates an extension of degree at most two over $\Bbb F_2\subset\bold f_{\goth p}$ when
the stabilizer $\Gamma_v=g(G_0)g^{-1}\cap\Gamma_0(\goth p)$ is non-trivial. Using the identity:
$$\left(\matrix0&1\\1&1+u\endmatrix\right)=
\left(\matrix0&1\\1&u\endmatrix\right)\cdot
\left(\matrix1&1\\0&1\endmatrix\right)$$
and the fact that $\left(\smallmatrix1&1\\0&1\endsmallmatrix\right)\in G_0$ we may assume that $u$ is either equal to $0$ or the element $t$ which we introduced in Definition 6.2. In the first case the vertex $v$ is equal to $v_0$ because $\left(\smallmatrix0&1\\ 1&0\endsmallmatrix\right)\in G_0$. In particular its stabilizer is non-trivial. In the second case the vertex $v$ is equal to $w$ by definition. Moreover there is a matrix $\left(\smallmatrix a&b\\ c&d\endsmallmatrix\right)\in GL_2(\Bbb F_2)$ such that $c\neq0$ and
$g\left(\smallmatrix a&b\\ c&d\endsmallmatrix\right)g^{-1}
\in\Gamma_0(\goth p)$. Therefore the stabilizer $\Gamma_w$ is non-trivial, too.\ $\square$
\enddefinition
For every vertex $v$ (or edge $e$) of $\Cal T$ let $\widetilde v$ (resp. $\widetilde e$) denote the image of $v$ (resp. $e$) in the factor graph $\Gamma_0(\goth p)\backslash\Cal T$. 
\proclaim{Lemma 6.6} Assume that $q=2$. Then the following holds:
\roster
\item"$(i)$" for every $n\in\Bbb N$ the degrees of the vertices
$\widetilde u_n$ and $\widetilde v_n$ are two,
\item"$(ii)$" if $\deg(\goth p)$ is even then the degree of
$\widetilde w$ is one.
\endroster
\endproclaim
\definition{Proof} Let $f$ be the unique generator of the ideal
$\goth p$. The matrix $w_{\goth p}=\left(\smallmatrix0&1\\ f&0\endsmallmatrix\right)\in GL_2(F)$ normalizes the subgroup $\Gamma_0(\goth p)$, hence its action on the Bruhat-Tits tree $\Cal T$ induces an automorphism of the factor graph $\Gamma_0(\goth p)\backslash\Cal T$. The square of $w_{\goth p}$ is a scalar matrix hence this automorphism is in fact an involution. As usual, we will call it the Atkin-Lehner involution. Note that:
$$\left(\matrix0&1\\ f&0\endmatrix\right)
\cdot\left(\matrix T^n&0\\ 0&1\endmatrix\right)=
\left(\matrix0&1\\ 1&0\endmatrix\right)\cdot
\left(\matrix T^{n+\deg(\goth p)}&0\\ 0&1\endmatrix\right)
\cdot\left(\matrix f\cdot T^{-\deg(\goth p)}&0\\ 0&1\endmatrix\right).$$
Since $\left(\smallmatrix f\cdot T^{-\deg(\goth p)}&0\\ 0&1\endsmallmatrix\right)\in GL_2(\Cal O_{\infty})$ we get that $\widetilde u_n$ is the image of $\widetilde v_n$ under the Atkin-Lehner involution. Hence it will be sufficient to verify claim $(i)$ only for the vertex $\widetilde v_n$. Because $3$ does not divide the order of the stabilizer $\Gamma_{v_n}$ the latter cannot act transitively on the set of edges with origin $v_n$. Hence the degree of $\widetilde v_n$ is at least two. For every $u\in\Bbb F_2$ we have:
$$\left(\matrix 1&T^nu\\ 0&1\endmatrix\right)
\cdot\left(\matrix T^n&0\\ 0&1\endmatrix\right)\cdot
\left(\matrix0&1\\ 1&0\endmatrix\right)=
\left(\matrix T^n&0\\ 0&1\endmatrix\right)\cdot
\left(\matrix u&1\\ 1&0\endmatrix\right).$$
By definition the product matrix in the equation above represents an edge on $\Cal T$ whose origin is $v_n$. Because $\left(\smallmatrix 1&T^nu\\ 0&1\endsmallmatrix\right)\in\Gamma_{v_n}$ it will be sufficient to show that these edges are different for different values of $u$. The latter follows from the fact that
$$\left(\matrix 0&1\\ 1&0\endmatrix\right)^{-1}
\cdot\left(\matrix 1&1\\ 1&0\endmatrix\right)=
\left(\matrix1&0\\ 1&1\endmatrix\right)\not\in Z(F_{\infty})
\Gamma_{\infty}.$$
The reduction of the element $t$ introduced in Definition 6.2 modulo
$\goth p$ satisfies the equation $x^2+x+1$. Hence the computations of the proof of Proposition 6.3 imply that the matrix $g=\left(\smallmatrix 0&1\\ 1&1\endsmallmatrix\right)\in\Gamma_{w}$. This matrix has order $3$. The group $G_0$ acts faithfully on the set of edges of $\Cal T$ with origin $v_0$. In particular every element of $G_0$ of order $3$ acts transitively on these edges. Therefore $g$ must act transitively on the set of edges of $\Cal T$ with origin
$w$.\ $\square$
\enddefinition
\definition{Notation 6.7} We will say that a half-infinite path $\gamma=\{e_0,\ldots,e_n,\ldots\}$ on an oriented graph $G$ is a half-line when for every vertex $v$ of $G$ there is at most one index $n\in\Bbb N$ such that $v=o(e_n)$. By a slight abuse of terminology we will call an oriented graph $G$ a half-line if there is a half-line $\gamma=\{e_0,\ldots,e_n,\ldots\}$ on $G$ such that for every $e\in\Cal E(G)$ either $e$ or $\overline e$ is listed in $\gamma$. Note that the graph $G$ determines uniquely the half-line $\gamma$. Let $H$ be a subgraph of the oriented graph $G$ which is a half-line and let $\{e_0,\ldots,e_n,\ldots\}$ be the half-line on $G$ which is uniquely determined by $H$. We will say that $H$ is a maximal half-line if $o(e_0)$ is the only vertex of $H$ whose degree in $G$ is not equal to two. 
\enddefinition
\proclaim{Lemma 6.8} The graph $\Gamma_0(\goth p)\backslash\Cal T$ decomposes uniquely as the edge-disjoint union of a finite graph and two maximal half-lines.
\endproclaim
\definition{Proof} By Theorem II.9 of [29] the graph $\Gamma_0(\goth p)\backslash\Cal T$ is the edge-disjoint union of a finite graph and finitely many half-lines. Moreover there is a bijection between the latter set of half-lines and the set $\Gamma_0(\goth p)\backslash\Bbb P^1(F)$ where $\Gamma_0(\goth p)$ acts on $\Bbb P^1(F)$ via the M\"obius action. The set $\Gamma_0(\goth p)\backslash\Bbb P^1(F)$ has two elements. Hence $\Gamma_0(\goth p)\backslash\Cal T$ is the edge-disjoint union of a finite graph and two half-lines. The graph
$\Gamma_0(\goth p)\backslash\Cal T$ has vertices whose degree is neither one nor two, for example because the group $H_!(\Cal T,\Bbb Z)^{\Gamma_0(\goth p)}$ is not trivial by Proposition 4.10. Therefore each of the two half-lines above is contained in a unique maximal half-line.\ $\square$
\enddefinition
Let $\Cal G(\goth p)$ denote the finite graph in the proposition above. We say that a vertex $v$ of $\Gamma_0(\goth p)\backslash\Cal T$ is cuspidal if it is a vertex of one of the two maximal half-lines of $\Gamma_0(\goth p)\backslash\Cal T$ but it is not a vertex of $\Cal G(\goth p)$. In the rest of this chapter we assume that $q=2$. 
\proclaim{Proposition 6.9}  Let $v$ be a vertex of $\Gamma_0(\goth p)\backslash\Cal T$. Then $v$ is cuspidal if and only if $v$ is equal to one of the vertices $\widetilde u_n$ or
$\widetilde v_n$ for some $n\in\Bbb N$.
\endproclaim
\definition{Proof} Every cuspidal vertex of $\Gamma_0(\goth p)\backslash\Cal T$ has degree two. In particular its stabilizer is non-trivial. Hence $v$ is equal to one of the vertices $\widetilde u_n$ or $\widetilde v_n$ for some $n\in\Bbb N$ if it is cuspidal by Proposition 6.3 and by claim $(i)$ of Lemma 6.6. On the other hand it will be sufficient to prove the converse for the vertices $\widetilde u_n$ only because the vertices $\widetilde v_n$ are images of these under the Atkin-Lehner involution. Let $e_n$ denote the edge of the Bruhat-Tits tree represented by the matrix:
$$\left(\matrix T&0\\ 0&T^{n+\deg(\goth p)}\endmatrix\right)=
\left(\matrix0&1\\ T^{n+\deg(\goth p)}&0\endmatrix\right)
\cdot\left(\matrix 0&1\\ T&0\endmatrix\right)$$
for every natural number $n\in\Bbb N$. By definition we have that $t(e_n)=u_n$ for every $n\in\Bbb N$. On the other hand:
$$\left(\matrix T&0\\ 0&T^{n+\deg(\goth p)}\endmatrix\right)=
\left(\matrix0&1\\ T^{n+\deg(\goth p)-1}&0\endmatrix\right)
\cdot\left(\matrix 0&T\\ T&0\endmatrix\right).$$
Because $\left(\smallmatrix0&T
\\ T&0\endsmallmatrix\right)\in Z(F_{\infty})GL_2(\Cal O_{\infty})$ 
we get that $o(e_n)=u_{n-1}$ for every $n\in\Bbb N$. Therefore the sequence $\{e_0,\ldots,e_n,\ldots\}$ is a half-infinite path. Because the vertices $o(\widetilde e_n)=\widetilde u_{n-1}$ are all different the sequence $\{\widetilde e_0,\ldots,\widetilde e_n,\ldots\}$ is a 
half-line. The claim is now clear.\ $\square$
\enddefinition
\definition{Notation 6.10} Note that the stabilizer of two vertices (or edges) of $\Cal T$ in $\Gamma_0(\goth p)$ which are equivalent under the action of this group are conjugate subgroups. In particular the isomorphism class of the stabilizer of a vertex $v$ (or and edge $e$) only depends on the image of $v$ (respectively $e$) in the quotient graph $\Gamma_0(\goth p)\backslash\Cal T$. Therefore we may talk about the stabilizer of a vertex (or an edge) of $\Gamma_0(\goth p)\backslash\Cal T$ being non-trivial (in the group $\Gamma_0(\goth p)$). 
\enddefinition
\proclaim{Proposition 6.11} Let $v$ be a vertex of $\Gamma_0(\goth p)\backslash\Cal T$ which is not cuspidal and whose stabilizer is non-trivial. Then the following holds:
\roster
\item"$(i)$" there is no vertex other than $v$ whose stabilizer is non-trivial
and which is not cuspidal,
\item"$(ii)$" the degree of $v$ is one,
\item"$(iii)$" the vertex $v$ is fixed by the Atkin-Lehner involution,
\item"$(iv)$" the unique neighbor of $v$ has trivial stabilizer.
\endroster
\endproclaim
\definition{Proof} The first claim is an immediate corollary of Propositions 6.3 and 6.9. The second claim follows from claim $(ii)$ of Lemma 6.6. By claims $(i)$ and $(ii)$ the vertex $v$ is the only
vertex of $\Gamma_0(\goth p)\backslash\Cal T$ whose degree is one. Hence it must be fixed by every automorphism of this graph. Therefore claim $(iii)$ is true. Assume that the unique neighbor $w'$ of $v$ has non-trivial stabilizer. Then this neighbor must be cuspidal. That means that the maximal half-line to which $w'$ belongs to terminates in a degree one vertex. In this case the graph $\Gamma_0(\goth p)\backslash\Cal T$ cannot be connected
which is a contradiction. Claim $(iv)$ in now clear.\ $\square$
\enddefinition
We will say that an edge $e$ of $\Gamma_0(\goth p)\backslash\Cal T$ is cuspidal if it is an edge of one of the two maximal half-lines of
$\Gamma_0(\goth p)\backslash\Cal T$. For every edge $e\in\Cal E(\Cal T)$ the stabilizer of $e$ in $\Gamma_0(\goth p)$ is the intersection of the stabilizers of the vertices $o(e)$ and $t(e)$. Hence claim $(iv)$ above has the following immediate
\proclaim{Corollary 6.12} Let $e$ be an edge of $\Gamma_0(\goth p)\backslash\Cal T$ which is not cuspidal. Then the stabilizer of $e$ is trivial.\ $\square$
\endproclaim
The next proposition will be used in the proof of the formal
immersion property.
\proclaim{Proposition 6.13} The following holds:
\roster
\item"$(i)$" the unique common vertex of the graph $\Cal G(\goth p)$ and the maximal half-line of $\Gamma_0(\goth p)\backslash\Cal T$ containing the vertices
$\{\widetilde u_n|n\in\Bbb N\}$ is $\widetilde u_{-1}$,
\item"$(ii)$" there is a harmonic cochain $\phi\in\Lambda(\goth p)$ whose reduction modulo $2$, considered as a function $\phi:
\Cal E(\Gamma_0(\goth p)\backslash\Cal T)\rightarrow\Bbb F_2$, is zero on exactly one edge with origin $\widetilde u_{-1}$. 
\endroster
\endproclaim
\definition{Proof} The unique common vertex is $o(\widetilde e_0)$ where $e_0$ is the edge introduced in the proof of Proposition 6.9. Claim $(i)$ is now clear. By Lemma 4.8 we know that:
$$(1+q-T_{(T)})[0,\infty]=[0,1/T]\in\Lambda(\goth p).$$
We claim that $\phi=[0,1/T]$ satisfies the condition in claim $(ii)$. Let $f_n$ denote the edge of the Bruhat-Tits tree represented by the matrix:
$$\left(\matrix0&T^n\\ 1&T^{n+1}\endmatrix\right)=
\left(\matrix0&1\\ 1&T\endmatrix\right)
\cdot\left(\matrix 1&0\\ 0&T^n\endmatrix\right)$$
for every integer $n\in\Bbb Z$. Note that:
$$\left(\matrix0&T^n\\ 1&T^{n+1}\endmatrix\right)
\cdot\left(\matrix 0&1\\ T&0\endmatrix\right)=
\left(\matrix 0&T^{n+1}\\ 1&T^{n+2}\endmatrix\right)
\cdot\left(\matrix 0&1\\ 1&0\endmatrix\right).$$
Because $\left(\smallmatrix0&1
\\ 1&0\endsmallmatrix\right)\in GL_2(\Cal O_{\infty})$ 
we get that $t(f_n)=o(f_{n+1})$ for every $n\in\Bbb Z$. Hence the sequence $\gamma=\{\ldots,f_{-n},\ldots,f_n,f_{n+1},\dots\}$ is an infinite path. It is clear that we have
$S(0,1/T)=\{f_n|n\in\Bbb Z\}$. Therefore $\gamma$ is the fundamental arch $\overline{0T^{-1}}$ by definition. Since
$$\left(\matrix 1&0\\ 0&T^n\endmatrix\right)=
\left(\matrix T^{-n}&0\\ 0&1\endmatrix\right)
\cdot\left(\matrix T^n&0\\ 0&T^n\endmatrix\right)$$
we get that $o(f_n)=\left(\smallmatrix0&1
\\ 1&T\endsmallmatrix\right)\cdot v_{-n}$ for every negative integer
$n$. Moreover:
$$\left(\matrix0&T^n\\ 1&T^{n+1}\endmatrix\right)=
\left(\matrix 1&0\\ T&1\endmatrix\right)
\cdot\left(\matrix T^n&0\\ 0&1\endmatrix\right)
\cdot\left(\matrix0&1\\ 1&0\endmatrix\right)$$
and $\left(\smallmatrix0&1
\\ 1&0\endsmallmatrix\right)\in GL_2(\Cal O_{\infty})$ hence we have
$o(f_n)=\left(\smallmatrix1&0\\ T&1\endsmallmatrix\right)\cdot
v_n$ for every non-negative integer $n$. 

For the sake of simple notation let $d$ denote $\deg(\goth p)$. From Lemmas 6.4 and 6.5 it follows that for every $g\in GL_2(A)-\Gamma_0(\goth p)$ and natural number $n\geq d-1$ the vertex $g\cdot v_n$ is equivalent to $u_n$ under the action of $\Gamma_0(\goth p)$. Therefore we have $o(\widetilde f_n)=\widetilde u_{-1}$ if and only if either $n=d-1$ or $n=1-d$, and we have $t(\widetilde f_n)=\widetilde u_{-1}$ if and only if  either $n=d-2$ or $n=-d$. Moreover both
$o(\widetilde f_{-d})$ and $t(\widetilde f_{d-1})$ are equal to
$\widetilde u_0$ hence the contribution of the edges
$\widetilde f_{-d}$ and $\widetilde f_{d-1}$ to $[0,1/T]$ cancel each other out. Now we only have to show that $t(\widetilde f_{1-d})$ and
$o(\widetilde f_{d-2})$ are different. Let $u\in A$ be the unique element such that $\deg(u)<d$ and $uT\equiv1\mod\goth p$. Since we have:
$$\left(\matrix 1&0\\ T&1\endmatrix\right)=
\left(\matrix T^{-1}&1\\ 0&T\endmatrix\right)
\cdot\left(\matrix 0&1\\ 1&T^{-1}\endmatrix\right)$$
the reductions of the matrices $\left(\smallmatrix 1&0\\ T&1\endsmallmatrix\right)$ and $\left(\smallmatrix0&1\\ 1&u\endsmallmatrix\right)$ modulo $\goth p$ lie in the same left $B(\bold f_{\goth p})$-coset. As we saw in the proof of Lemma 6.4 the reductions of the elements of $M(\goth p)$ form a complete system of representatives for the left cosets of $B(\bold f_{\goth p})$ hence the matrices $\left(\smallmatrix 1&0\\ T&1\endsmallmatrix\right)$ and $\left(\smallmatrix 0&1\\ 1&u\endsmallmatrix\right)$ lie in the same left $\Gamma_0(\goth p)$-coset. Therefore $o(f_{d-2})$ is equivalent to
$\left(\smallmatrix0&1\\ 1&u\endsmallmatrix\right)\cdot v_{d-2}$ modulo the action of $\Gamma_0(\goth p)$. By Lemma 6.5 the claim now follows from the fact that $\deg(u)=d-1$.\ $\square$
\enddefinition
Let $h_0$ and $h_{\infty}$ denote the half-lines $\{e_0,\ldots,e_n,\ldots\}$ and $\{w_{\goth p}\cdot e_0,\ldots,w_{\goth p}\cdot e_n,\ldots\}$ on $\Cal T$, respectively, where we continue to use the notation of the proofs of Lemma 6.6 and Proposition 6.9. For $i=0$, $\infty$ let $H_i$ denote the subgraph of $\Cal T$ which is the half-line determined by $h_i$. For every $i\in\Bbb P^1(F)$ let $\text{Stab}_i$ denote the subgroup of those elements of $\Gamma_0(\goth p)$  which leave $i$ fixed where we let $\Gamma_0(\goth p)$ act on $\Bbb P^1(F)$ via the usual M\"obius action. 
\proclaim{Lemma 6.14} Let $i$ be an index such that either $i=0$ or $i=\infty$ and let $\gamma$ be an element of $\Gamma_0(\goth p)$. Then the following holds:
\roster
\item"$(i)$" we have $\gamma(H_0)\cap H_{\infty}=\emptyset$,
\item"$(ii)$" if $\gamma(H_i)\cap H_i\neq\emptyset$ then $\gamma\in\text{\rm Stab}_i$.
\endroster
\endproclaim
\definition{Proof} The images of $H_0$ and $H_{\infty}$ in the quotient graph $\Gamma_0(\goth p)$ are disjoint maximal half-lines, hence claim $(i)$ is clear. Assume now that $\gamma(H_i)\cap H_i\neq\emptyset$ for some $i$. Because the matrix $w_{\goth p}$ interchanges the half-lines $h_0$ and $h_{\infty}$ as well as the points $0$ and $\infty$ on $\Bbb P^1(F)$ via its M\"obius action we may assume that $i=0$ by conjugating $\gamma$ by $w_{\goth p}$, if necessary. Let $n\in\Bbb N\cup\{-1\}$ be an index such that $u_n$ lies in $\gamma(H_i)\cap H_i$. Then $u_n=\gamma(u_m)$ for some $m\in\Bbb N\cup\{-1\}$. Because $\widetilde u_n$ is different from $\widetilde u_m$ for every $m\neq n$ we must have $u_n=\gamma(u_n)$. Therefore we have:
$$\gamma\in\left(\matrix0&1\\ 1&0\endmatrix\right)
G_{n+d}\left(\matrix0&1\\ 1&0\endmatrix\right)^{-1}
\!\!\!\!\cap\Gamma_0(\goth p)
\subset
\{\left(\matrix1&0\\ u&1\endmatrix\right)|u\in\goth p\}
=\text{Stab}_0\text{.\ $\square$}$$
\enddefinition

\heading 7. The regular model $\goth X_0(\goth p)$
\endheading

\definition{Definition 7.1} Let $\Omega$ denote the rigid analytic upper half plane, or Drinfeld's upper half plane over $F_{\infty}$. The set of ($\Bbb C_{\infty}$-valued) points of $\Omega$ is $\Bbb C_{\infty}-F_{\infty}$, denoted also by $\Omega$ by abuse of notation. The group $GL_2(F_{\infty})$ acts on Drinfeld's upper half plane on the left via M\"obius transformations. The restriction of this action to $\Gamma_0(\goth n)$ is discrete, hence the set $\Gamma_0(\goth n)\backslash\Omega$ has naturally the structure of a rigid analytic curve. Let $Y_0(\goth n)$ also denote the underlying rigid analytical space of the base change of $Y_0(\goth n)$ to $F_{\infty}$ by abuse of notation.
\enddefinition
\proclaim{Theorem 7.2} There is a rigid-analytical isomorphism:
$$Y_0(\goth n)\cong\Gamma_0(\goth n)\backslash\Omega.$$
\endproclaim
\definition{Proof} See [2], Theorem 6.6.\ $\square$
\enddefinition
\definition{Definition 7.3} Let $X$ be a rigid analytic space over $F_{\infty}$ and assume that $\Cal U$ is an admissible covering of $X$ by affinoids which is pure in the sense of [11] (for the definition of purity and the following concepts see 2.10 of [11] on page 116). We may associate a formal scheme $\goth X$ to the data $(X,\Cal U)$ whose generic fiber in the sense of Raynaud is canonically isomorphic to $X$. The topological space underlying $\goth X$ is the same as the topological space underlying the reduction $\overline X_{\Cal U}$ of the rigid analytic space associated to $X$ with respect to the cover $\Cal U$. For every element $U$ of the cover $\Cal U$ let $\overline U$ denote its reduction in the sense used in [11]. By definition both $\overline U$ and $\overline X_{\Cal U}$ are schemes over $\bold f_{\infty}$ and 
the natural map $\overline U\rightarrow\overline X_{\Cal U}$ is an open immersion. Moreover for every $U$ as above let $\goth U$ denote the formal scheme:
$$\goth U=Spf(\{f\in\Cal O(U)|\|f\|\leq1\}),$$
where $Spf(\cdot)$ denotes the formal spectrum and $\|\cdot\|$ denotes the supremum norm on $\Cal O(U)$ with respect to some absolute value $|\cdot|$ on $\Bbb C_{\infty}$ corresponding to the valuation on this field. Then the topological spaces underlying $\goth U$ and $\overline U$ are canonically homeomorphic. By definition the formal scheme $\goth X$ is just the unique locally ringed space whose underlying topological space is $\overline X_{\Cal U}$ such that the restriction of its structure sheaf $\Cal O_{\goth X}$ to the open $\overline U$ is the structure sheaf $\Cal O_{\goth U}$ of the formal scheme $\goth U$ for every $U$ as above.
\enddefinition
\definition{Definition 7.4} Let $U(1)$ denote the rational subdomain of $\Omega$ whose underlying set of points is:
$$U(1)=\{z\in\Bbb C_{\infty}||\pi|\leq|z|\leq1,\ |z-c|\geq1,|z-c\pi|\geq|\pi|\
(\forall c\in\bold f^*_{\infty})\},$$
where $\pi\in F_{\infty}$ is a uniformizer. For every $\rho\in GL_2(F_{\infty})$ and $z\in\Bbb P^1$ let $\rho(z)$ denote the image of $z$ under the M\"obius transformation corresponding to $\rho$. For every $\rho$ as above let $U(\rho)$ denote the rational subdomain whose underlying set of points is:
$$U(\rho)=\{z\in\Bbb C_{\infty}|\rho^{-1}(z)\in U(1)\}.$$
The system $\Cal U=\{U(\rho)|\rho\in GL_2(F_{\infty})\}$ is a pure covering of $\Omega$ by affinoid subdomains. Let $\widehat{\Omega}$ denote the corresponding formal scheme. Because the cover $\Cal U$ is invariant with respect to the M\"obius action the formal scheme
$\widehat{\Omega}$ is equipped with an action of $GL_2(F_{\infty})$ such that for every $\rho\in GL_2(F_{\infty})$ the automorphism of $\Omega$ induced by the action of $\rho$ on $\widehat{\Omega}$
is the usual M\"obius action of $\rho$.
\enddefinition
\definition{Definition 7.5} By slightly extending the usual terminology we will call a scheme $C$ defined over a field a curve if it is reduced,
locally of finite type and of dimension one. A curve $C$ is said to
have normal crossings if every singular point of $C$ is an ordinary
double point in the usual sense. For any curve $C$ with normal
crossings let $\widetilde C$ denote its normalization, and let
$\widetilde S(C)$ denote the pre-image of the set $S(C)$ of singular
points of $C$. The incidence graph of $C$ is the oriented graph whose set of vertices is the set of irreducible components of $\widetilde C$,
and its set of edges is the set $\widetilde S(C)$ such that the
original vertex of any edge $x\in\widetilde S(C)$ is the irreducible
component of $\widetilde C$ which contains $x$ and the terminal
vertex of $x$ is the irreducible component which contains the unique
other element $\overline x$ of $\widetilde S(C)$ which maps with
respect to the normalization map to the same
 singular point as $x$.
\enddefinition
\definition{Remark 7.6} The special fiber $\widehat{\Omega}_{\bold f_{\infty}}$ of the formal scheme $\widehat{\Omega}$ is a curve with normal crossings over $\bold f_{\infty}$. The incidence graph of $\widehat{\Omega}_{\bold f_{\infty}}$ is canonically isomorphic to the Bruhat-Tits tree $\Cal T$ in the sense that the natural action of $GL_2(F_{\infty})$ on the incidence graph is the usual action of $GL_2(F_{\infty})$ on the Bruhat-Tits tree under this identification. (We will give a more explicit description of this bijection in Definition 7.8 below.) We may form the quotient formal scheme $\Gamma_0(\goth n)\backslash\widehat{\Omega}$. By the above the incidence graph of the special fiber of the latter is $\Gamma_0(\goth n)\backslash\Cal T$, which is not finite, hence the reduction of the formal scheme $\Gamma_0(\goth n)\backslash\widehat{\Omega}$ is not even of finite type. In particular this formal scheme cannot be the formal completion of an algebraic curve over Spec$(\Cal O_{\infty})$ along its special fiber. Nevertheless there is a way to obtain a model from this formal scheme, using the work of Reversat in [23]. Next we will give our own account of this theory. Another possible reference is section 4.2 of [19].
\enddefinition
\definition{Notation 7.7} Assume now that $q=2$. For the sake of simple notation let $\Gamma$ denote $\Gamma_0(\goth p)$ in the rest of this chapter. Let $\Gamma_{tor}$ be the normal subgroup of $\Gamma$ generated by the torsion elements, and let $\Gamma'=\Gamma/\Gamma_{tor}$. For $i=0$, $\infty$ let $H_i'$ denote the subgraph of $\Cal T$ we get from $H_i$ by removing its only vertex of degree one, where $H_i$ is the graph which was introduced before Lemma 6.14. Moreover let $G_i$ and $G_i'$ denote the subgraph of $\Cal T$ which is the union of the graphs $\gamma(H_i)$ and $\gamma(H_i')$, respectively, where $\gamma$ is an arbitrary element of $\text{Stab}_i$. If $\gamma$ is an element of $\Gamma$ such that $\gamma(G_i)\cap G_i\neq\emptyset$ then $\gamma\in\text{Stab}_i$ by Lemma 6.14. Similarly $\gamma(G_0)\cap G_{\infty}=\emptyset$ for every $\gamma\in\Gamma$.
\enddefinition
\definition{Definition 7.8} By slight abuse of notation for every ringed space such as $\Omega$ or $\widehat\Omega_{\bold f_{\infty}}$ let the same symbol denote the underlying topological space, too. Let red:$\Omega\rightarrow\widehat\Omega_{\bold f_{\infty}}$ be the reduction map. It is a continuous $GL_2(F_{\infty})$-equivariant map between the underlying topological spaces. For every irreducible component $v$ of $\widehat\Omega_{\bold f_{\infty}}$ let $ns(v)$ denote the topological subspace which we get from $v$ by removing its singular points. For every edge $e$ of the incidence graph of $\widehat\Omega_{\bold f_{\infty}}$ let $\widehat e$ denote the image of $e$ under the normalization map. For every subgraph $G$ of the incidence graph of $\widehat\Omega_{\bold f_{\infty}}$ the realization of $G$, denoted by $r(G)$, is the subspace of the topological space
$\widehat\Omega_{\bold f_{\infty}}$ which is the union of the subspaces $ns(v)$, where $v$ is any element of $\Cal V(G)$, and the points $\widehat e$, where $e$ is any element of $\Cal E(G)$. For every edge $e$ of the incidence graph of
$\widehat\Omega_{\bold f_{\infty}}$ let $G(e)$ denote the subgraph whose vertices are $o(e)$ and $t(e)$ and whose edges are $e$ and
$\overline e$. For every $\rho\in\ GL_2(F_{\infty})$ let $e(\rho)$
denote the edge of $\Cal T$ represented by $\rho$. Under the identification between $\Cal T$ and the incidence graph of
$\widehat\Omega_{\bold f_{\infty}}$ which we mentioned in Remark 7.6 the pre-image of $r(G(e(\rho)))$ under the reduction map is the domain $U(\rho)$ for every $\rho\in\ GL_2(F_{\infty})$.
\enddefinition
\definition{Construction 7.9} For $i=0$, $\infty$ let $\Omega_i$ and
$\Omega_i'$ be the pre-image of $r(H_i)$ and $r(H_i')$, respectively, under the reduction map. By definition $\Omega_i'$ is an open subdomain of $\Omega_i$ and these domains are invariant under the action of $\text{Stab}_i$. Moreover it is clear from the remark at the end of Notation 7.7 that if we have $\gamma(\Omega_i)\cap\Omega_j\neq\emptyset$ for some $\gamma\in\Gamma$ then $i=j$ and $\gamma\in\text{Stab}_i$. Hence the natural map $j_i:\text{Stab}_i\backslash\Omega_i\rightarrow \Gamma_{tor}\backslash\Omega$ is an open immersion. By Lemma 2.4 of [23] on page 384 for every $i$ the quotient $\text{Stab}_i\backslash\Omega_i'$ is analytically isomorphic to a unit disc $D_i$ with the origin removed. Let $B_i$ be the rigid analytic space $\text{Stab}_i\backslash\Omega_i\cup D_i$ obtained by gluing $D_i$ with $\text{Stab}_i\backslash\Omega_i$ along their common admissible subspace $\text{Stab}_i\backslash\Omega_i'$. For every $\gamma\in\Gamma'$ let $B_i^{\gamma}$ be a copy of
$B_i$. Denote by $\Sigma$ the rigid analytic space
$$\Sigma=(\Gamma_{tor}\backslash\Omega)\cup
\bigcup_{\gamma\in\Gamma'}B^{\gamma}_0\cup
\bigcup_{\gamma\in\Gamma'}B^{\gamma}_{\infty}$$
obtained by gluing each $B^{\gamma}_i$ with $\Gamma_{tor}\backslash\Omega$ along their common admissible subspace
$\gamma(j_i(\text{Stab}_i\backslash\Omega_i))$. Let $f:\Omega
\rightarrow\Gamma_{tor}\backslash\Omega\subset\Sigma$ be the canonical morphism of analytic spaces. By Lemma 2.6 of [23] on page 385 the system:
$$\Cal V=\{f(U(\rho))|e(\rho)\not\in\Cal E(G_0)\cup\Cal E(G_{\infty})\}
\cup\{B^{\gamma}_0|\gamma\in\Gamma'\}
\cup\{B^{\gamma}_{\infty}|\gamma\in\Gamma'\}$$
is a pure admissible covering of $\Sigma$ by affinoids. Let $\widehat\Sigma$ denote the formal scheme associated to the data $(\Sigma,\Cal V)$. The natural action of the group $\Gamma'$ on $\Gamma_{tor}\backslash\Omega$ extends uniquely to a discontinuous action on $\Sigma$. The group $\Gamma'$ leaves the cover $\Cal V$ invariant hence we may equip $\widehat\Sigma$ with an action of
$\Gamma'$ such that the induced action on $\Sigma$ is the one above. 
\enddefinition
By construction the special fiber of $\Gamma\backslash\widehat\Sigma$ has incidence graph $\Cal G(\goth p)$. In particular this graph is finite so the special fiber of $\Gamma\backslash\widehat\Sigma$ is proper. Hence by Grothendieck's algebraization theorem the formal scheme $\Gamma\backslash\widehat\Sigma$ is the formal completion of a projective curve $\goth X_0(\goth p)\rightarrow\text{Spec}(\Cal O_{\infty})$ along its special fiber.
According to the main result of [23] (Theorem 2.7 on page 385):
\proclaim{Proposition 7.10} The curve $\goth X_0(\goth p)$ is a totally degenerate semi-stable model of $X_0(\goth p)$ over the spectrum of $\Cal O_{\infty}$.\ $\square$
\endproclaim
\definition{Notation 7.11} By slight abuse of notation let the index $i$ also denote the unique point of $B_i$ which does not lie in $\text{Stab}_i\backslash\Omega_i$. Under the action of $\Gamma'$ the corresponding points in the copies $B^{\gamma}_i$ are all equivalent. Let $i$ denote the common image of these points in the quotient
$\Gamma'\backslash\Sigma$ as well. The identification of the rigid analytical space underlying the base change of $X_0(\goth n)$ to $F_{\infty}$ and the quotient $\Gamma'\backslash\Sigma$ identifies the two cusps of $X_0(\goth p)$ and the points $0$ and $\infty$. We will use this identification without further notice.
\enddefinition
\proclaim{Proposition 7.12} The model $\goth X_0(\goth p)$ is regular.
\endproclaim
\definition{Proof} Because $\goth X_0(\goth p)$ is semi-stable and its base change to Spec$(F_{\infty})$ is smooth, it is regular in the complement of the ordinary double points in the special fiber. The system:
$$\Cal V_0=\{f(U(\rho))|e(\rho)\not\in\Cal E(G_0)\cup\Cal E(G_{\infty})\}$$
is a pure admissible covering by affinoids of the domain $\Sigma_0\subset\Sigma$ it covers. Let $\widehat\Sigma_0$ denote the open formal subscheme of $\widehat\Sigma$ corresponding to the data
$(\Sigma_0,\Cal V_0)$. By construction $\Sigma_0$ is invariant under the action of $\Gamma'$ and the open formal subscheme $\Gamma'
\backslash\widehat\Sigma_0$ of $\goth X_0(\goth p)$ contains every ordinary double point in the special fiber. The formal scheme
$\Gamma'\backslash\widehat\Sigma_0$ is isomorphic to the quotient $\Gamma\backslash\widehat\Omega_0$ where $\widehat\Omega_0$ is the unique $\Gamma$-invariant open formal subscheme of
$\widehat\Omega$ whose generic fiber in the sense of Raynaud is the domain $\Omega_0\subset\Omega$ covered by the system:
$$\Cal U_0=\{U(\rho)|
e(\rho)\not\in\bigcup_{\gamma\in\Gamma}\Cal E(\gamma(G_0))\cup\bigcup_{\gamma\in\Gamma}\Cal E(\gamma(G_{\infty}))\}.$$
By Corollary 6.12 the stabilizer of every edge of $\Cal G(\goth p)$ is trivial. Hence the formal completion of the local ring at every ordinary double point $v$ of the special fiber of $\Gamma\backslash\widehat\Omega_0$ is isomorphic to the formal completion of the local ring at any ordinary double point of the special fiber of $\widehat\Omega$ which maps to $v$ with respect to the covering $\widehat\Omega\rightarrow\Gamma\backslash\widehat\Omega$ by the above. This completion is isomorphic to the ring $\Cal O_{\infty}[[x,y]]/(xy-\pi)$, where $\pi$ is again a uniformizer of $F_{\infty}$. As the latter is regular the claim is now clear.\ $\square$
\enddefinition
\definition{Notation 7.13} Because the scheme $\goth X_0(\goth p)$ over Spec$(\Cal O_{\infty})$ is projective, every $F_{\infty}$-valued point $P:\text{Spec}(F_{\infty})\rightarrow X_0(\goth p)$ extends uniquely to an $\Cal O_{\infty}$-valued point $\text{Spec}(\Cal O_{\infty})\rightarrow\goth X_0(\goth p)$ which we will denote by the same symbol by the usual abuse of notation. For every point $P\in X_0(\goth p)$ we define the reduction of $P$ as the closed point which is the intersection of the section $P:\text{Spec}(\Cal O_{\infty})\rightarrow\goth X_0(\goth p)$ and the special fiber of $\goth X_0(\goth p)$. Let $\Cal J_0(\goth p)$ and $\Cal W(\goth p)$ denote the N\'eron model over $\Bbb P^1_{\Bbb F_q}$ of the Jacobian $J_0(\goth p)$ and the winding quotient $W(\goth p)$, respectively. By slight abuse of notation we are also going to let the same symbol denote the base change of this model to Spec$(\Cal O_v)$, which is also the N\'eron model over this base, for every place $v$ of $F$. Let $\goth X_0^*(\goth p)$ denote the Spec$(\Cal O_{\infty})$-scheme we get from $\goth X_0(\goth p)$ by removing the ordinary double points in the special fiber. Because $\goth X_0^*(\goth p)$ is smooth over Spec$(\Cal O_{\infty})$ there is a map $\goth X^*_0(\goth p)\rightarrow\Cal J_0(\goth p)$ whose base change to Spec$(F_{\infty})$ is the Albanese imbedding $X_0(\goth p)\hookrightarrow J_0(\goth p)$ with base point $0$. Let $\iota:\goth X^*_0(\goth p)\rightarrow\Cal W(\goth p)$ denote the composition of this map and the morphism $h:\Cal J_0(\goth p)\rightarrow\Cal W(\goth p)$ induced by the quotient map $J_0(\goth p)\rightarrow W(\goth p)$.
\enddefinition
\proclaim{Proposition 7.14} The map $\iota:\goth X^*_0(\goth p)\rightarrow\Cal W(\goth p)$ is a formal immersion at $0$.
\endproclaim
\definition{Proof}  For every scheme $S$ over Spec$(\Cal O_{\infty})$ let $S_{\infty}$ denote its base change to Spec$(\bold f_{\infty})$. For any algebraic group $G$ defined over a field $K$ let Cot$(G)$ denote the $K$-linear cotangent space of $G$ at the identity element. Grothendieck duality furnishes an isomorphism:
$$\Theta:\text{Cot}(\Cal J_0(\goth p)_{\infty})
\rightarrow H^0(\goth X_0(\goth p)_{\infty},\omega)$$
where $\omega$ denotes the relative dualizing sheaf of the curve
$\goth X_0(\goth p)_{\infty}$ (see section $e$ of [14] on page 140). According to  section $b$ of [1] on page 77 the sheaf $\omega$ is the sheaf of $1$-forms $\eta$ on the normalization of $\goth X_0(\goth p)_{\infty}$ which are regular except for simple poles at the points in
$\widetilde S(\goth X_0(\goth p)_{\infty})$, and for every element $e$ of $\widetilde S(\goth X_0(\goth p)_{\infty})$ we have
Res${}_e(\eta)+\text{Res}_{\overline e}(\eta)=0$. By the residue theorem for every $\eta\in H^0(\goth X_0(\goth p)_{\infty},\omega)$ the map Res$(\eta):\widetilde S(\goth X_0(\goth p)_{\infty})\rightarrow\bold f_{\infty}$ given by the rule $e\mapsto\text{Res}_e(\eta)$ is a harmonic cochain on the incidence graph $\Cal G(\goth p)$ of $\goth X_0(\goth p)_{\infty}$. Because $\goth X_0(\goth p)_{\infty}$ is totally degenerate the homomorphism:
$$\text{Res}:H^0(\goth X_0(\goth p)_{\infty},\omega)
\rightarrow H(\Cal G(\goth p),\bold f_{\infty})$$
given by the rule $\eta\mapsto\text{Res}(\eta)$ is an isomorphism. 

Let $T$ be a split torus over the field $K$. Then Cot$(T)$ is canonically isomorphic to Hom$(T,\Bbb G_m)\otimes K$ where Hom$(T,\Bbb G_m)$ is the character group of $T$. Let $A$ be an abelian variety over $F_{\infty}$ which has a rigid-analytical uniformization by a split torus $U$ over $F_{\infty}$. Then the connected component of the identity of the special fiber of the N\'eron model of $A$ over Spec$(\Cal O_{\infty})$ is canonically isomorphic to the split torus $U'$ over $\bold f_{\infty}$ whose group of characters is Hom$(U,\Bbb G_m)$. Hence the cotangent space of the group variety $\Cal J_0(\goth p)_{\infty}$ and $\Cal W(\goth p)_{\infty}$ is isomorphic to
$\overline\Gamma_0(\goth p)\otimes\bold f_{\infty}$ and $\Lambda(\goth p)\otimes\bold f_{\infty}$, respectively, and the homomorphism $h^*:\text{Cot}(\Cal W(\goth p)_{\infty})
\rightarrow\text{Cot}(\Cal J_0(\goth p)_{\infty})$ induced by the morphism $h$ introduced in Notation 7.13 above is $\lambda\otimes\text{id}_{\bold f_{\infty}}$ under this identification where $\lambda:\Lambda(\goth p)\rightarrow\overline\Gamma_0(\goth p)$ is the inclusion map. 

For every $\eta\in H(\Cal G(\goth p),\bold f_{\infty})$ let the same symbol denote the extension of $\eta$ to a function $\Cal E(\Gamma_0(\goth p)\backslash\Cal T)\rightarrow\bold f_{\infty}$ which is zero on every edge not in $\Cal E(\Cal G(\goth p))$. This extension is a harmonic cochain on $\Gamma_0(\goth p)\backslash\Cal T$. The pull-back of every $\eta$ as above with respect to the quotient map
$\Cal E(\Cal T)\rightarrow\Cal E(\Gamma_0(\goth p)\backslash\Cal T)$ is a harmonic cochain by Corollary 6.12. Hence we may consider $H(\Cal G(\goth p),\bold f_{\infty})$ as a subgroup of $H(\Cal T,\bold f_{\infty})^{\Gamma_0(\goth p)}$. Under this identification the map:
$$\text{Res}\circ\Theta:\overline\Gamma_0(\goth p)\otimes\bold f_{\infty}\cong\text{Cot}(\Cal J_0(\goth p)_{\infty})
\rightarrow H(\Cal G(\goth p),\bold f_{\infty})\subseteq
H(\Cal T,\bold f_{\infty})^{\Gamma_0(\goth p)}$$
is $R\otimes\text{id}_{\bold f_{\infty}}$ where $R$ denotes the reduction modulo $2$. According to the differential condition for being a formal immersion we need to show that there is an $\eta\in H^0(\goth X_0(\goth p)_{\infty},\omega)$ in the image of the map $\Theta\circ h^*$ which does not vanish at the reduction of $0$. The latter lies on the irreducible component $\widetilde u_{-1}$. Every non-zero $1$-form on the projective line $\Bbb P^1_{\bold f_{\infty}}$ which is regular except for simple poles at two closed points is actually non-zero at every point where it is regular. By the above the image of $\Theta\circ h^*$, considered as a subgroup of $H(\Cal T,\bold f_{\infty})^{\Gamma_0(\goth p)}$, is the reduction of $\Lambda(\goth p)$ modulo $2$. The proposition now follows from claim $(ii)$ of Proposition 6.13.\ $\square$
\enddefinition

\heading 8. Mod $\goth p$ Galois representations of rank two Drinfeld modules with good reduction everywhere
\endheading

\definition{Definition 8.1} At the beginning of this chapter we will not need to assume that $q=2$. For every Drinfeld module $\phi:A\rightarrow B\{\tau\}$, where $B$ is an $A$-algebra, and for every non-zero ideal $\goth n\triangleleft A$ let $\phi[\goth n]$ denote the $\goth n$-torsion group scheme of $\phi$ as usual. Recall that the latter is a finite, flat subgroup-scheme of $\Bbb G_a$ over Spec$(B)$. We say that an $A$-algebra $B$ has characteristic $\goth p$ if the annihilator of the $A$-module $B$ contains $\goth p$. This assumption implies that $B$ is an $\bold f_{\goth p}$-algebra. Similarly we will say that a Drinfeld module $\phi:A\rightarrow B\{\tau\}$ has characteristic $\goth p$ if $B$ has characteristic $\goth p$. Let $f\in\goth p$ be a generator of the ideal $\goth p$.
\enddefinition
\proclaim{Lemma 8.2} For every field extension $\bold k$ of $\bold f_{\goth p}$ and for every rank two Drinfeld module $\phi: A\rightarrow\bold k\{\tau\}$ of characteristic $\goth p$ the degree of the lowest non-zero term of the polynomial:
$$\phi(f)=\sum_{n=0}^{2\deg(\goth p)}a_n\tau^n\in\bold k\{\tau\}$$
is either $\deg(\goth p)$ or $2\deg(\goth p)$.
\endproclaim
\definition{Proof} According to the remark following Proposition 5.1 of [3], page 178 the polynomial $\phi(f)$ has no terms of degree less than $\deg(\goth p)$. Recall that the subgroup-scheme $\phi[\goth p]$ of the additive group scheme over $\bold k$ is defined by the equation
$$\sum_{n=0}^{2\deg(\goth p)}a_nx^{q^{n\deg(\goth p)}}=0.$$
If the degree of the lowest non-zero term of $\phi(f)$ is not $\deg(\goth p)$ then the group scheme $\phi[\goth p]$ is connected by Satz 5.3 of [3], page 179. Hence $\phi(f)$ has no non-zero terms of degree less than $2\deg(\goth p)$. As the leading coefficient of $\phi(f)$ is the same as the leading coefficient of $f$ the claim is
now clear.\ $\square$
\enddefinition
\definition{Definition 8.3} Let $\phi:A\rightarrow\Cal O_{\goth p}\{\tau\}$ be a Drinfeld module of rank two over $\Cal O_{\goth p}$. Let $f\in\goth p$ be a generator of the ideal $\goth p$ as above and let
$$\phi(f)=\sum_{n=0}^{2\deg(\goth p)}a_n\tau^n\in\Cal O_{\goth p}\{\tau\}.$$
According to Lemma 8.2 we have:
$$\goth p(a_i)>0\quad\text{for }i<d,
\quad\text{and}\quad\goth p(a_d)=0,$$
where $d$ is either $\deg(\goth p)$ or $2\deg(\goth p)$, when we will say that $\phi$ has ordinary or supersingular reduction, respectively, following the usual terminology. For simplicity let $U$ denote the Gal$(\overline F_{\goth p}|F_{\goth p})$-module $\phi[\goth p](\overline F_{\goth p})$. Recall that the latter has the structure of a vector space over $\bold f_{\goth p}$ of dimension $2$ and the action of Gal$(\overline F_{\goth p}|F_{\goth p})$ is $\bold f_{\goth p}$-linear.
\enddefinition
\proclaim{Lemma 8.4} The Galois module $U$ is irreducible when
$\phi$ has supersingular reduction.
\endproclaim
\definition{Proof} Note that $t\in\overline F_{\goth p}$ is a non-zero element of $U$ if and only if it is a root of the polynomial:
$$g_f(x)=a_0+a_1x^{q-1}+\cdots+a_kx^{q^k-1}+\cdots+
a_{2\deg(\goth p)}x^{q^{2\deg(\goth p)}-1}.$$
Therefore it will be sufficient to prove that $g_f$ is irreducible in $F_{\goth p}[x]$. Because $a_0=f$ the latter follows at once from the
Sch\"onemann-Eisenstein criterion.\ $\square$
\enddefinition
\proclaim{Proposition 8.5} Assume that $\phi$ has ordinary reduction. Then the following holds:
\roster
\item"$(i)$" there is an exact sequence of Galois modules:
$$0@>>>U_0@>>>U@>>>U_1@>>>0$$
such that the Galois module $U_0$ is a one-dimensional
$\bold f_{\goth p}$-linear sub-vector space,
\item"$(ii)$" the Galois module $U_1$ is unramified. 
\endroster
\endproclaim
\definition{Proof} Let $\goth p:\overline F_{\goth p}^*\rightarrow\Bbb Q$ denote the extension of the valuation $\goth p$ to the separable closure $\overline F_{\goth p}$, too. Because $\phi$ has good reduction we have $\goth p(t)\geq0$ for every $t\in U$. Let $U_0$ denote the subset:
$$\{t\in U|\goth p(t)>0\}\subseteq U.$$
Clearly $0\in U_0$, moreover $U_0$ is closed under addition and left invariant
by the action of the absolute Galois group of $F_{\goth p}$. Let $g\in A$ be arbitrary. Then we have:
$$\phi(g)=\sum_{n=0}^{2\deg(g)}b_n\tau^n$$
for some $b_n\in\Cal O_{\goth p}$. Hence for every $t\in U_0$ we have $\goth p(\phi(g)(t))>0$ therefore $U_0$ is an $A$-module under the action of $A$ on $\overline F_{\goth p}$ induced by $\phi$. In particular $U_0$ is an $\bold f_{\goth p}$-subspace of $U$.
Let $U_1$ denote the quotient of $U$ by this Galois-invariant subspace. As we already noted $t\in\overline F_{\goth p}^*$ is an element of $U$ if and only if it is the root of the polynomial:
$$g_f(x)=a_0+a_1t^{q-1}+\cdots+a_kt^{q^k-1}+\cdots+
a_{2\deg(\goth p)}t^{q^{2\deg(\goth p)}-1}.$$
Let $d$ denote $\deg(\goth p)$ for simplicity. Because $\phi$ has ordinary reduction we know that
$$\goth p(a_0)=1,\ \goth p(a_kt^{q^k-1})\geq1+(q^k-1)\goth p(t)\ (k<d),
\ \goth p(a_dt^{q^d-1})=(q^d-1)\goth p(t),\text{ and}$$
$$\goth p(a_kt^{q^k-1})\geq(q^k-1)\goth p(t)\ (d<k\leq 2d)),$$
hence by comparing valuations we get that $\goth p(t)=(q^{\deg(\goth p)}-1)^{-1}$, if $t$ is an element of $U_0$, and $\goth p(t)=0$, otherwise. Because
$$f=a_0=\prod_{0\neq t\in U}t,$$
we get that the cardinality of $U_0$ is $q^{\deg(\goth p)}$ by comparing the valuations of both sides. Claim $(i)$ is now clear. Let
$\goth M$ and $\goth M^+$ denote the groups
$$\goth M=\{x\in\overline F_{\goth p}|\goth p(x)\geq0\}
\quad\text{and}
\quad\goth M^+=\{x\in\overline F_{\goth p}|\goth p(x)>0\},$$
respectively. The natural action of $\text{\rm Gal}(\overline F_{\goth p}|F_{\goth p})$ on  the quotient group $V_0=\goth M/\goth M^+$ is unramified by Proposition 7 of [27] on page 268. From the discussion above it is clear that the group $U_1$ injects into $V_0$. Claim $(ii)$ is now clear.\ $\square$
\enddefinition 
\definition{Definition 8.6} Recall that in Definition 4.2.1 of [12] on page 65 an $A$-submodule $\Lambda\subset\overline F_{\infty}$ is called an $F$-lattice if it satisfies the following properties:
\roster
\item"$(i)$" it is finitely generated and free as an $A$-module,
\item"$(ii)$" it is left invariant by the action of the absolute Galois group Gal$(\overline F_{\infty}|F_{\infty})$,
\item"$(iii)$" it is discrete, that is, its intersection with every open disk is finite.
\endroster
\enddefinition
\proclaim{Lemma 8.7} Let $\Lambda\subset\overline F_{\infty}$ be an $F$-lattice. Then the absolute Galois group $\text{\rm Gal}(\overline F_{\infty}|F_{\infty})$ of $F_{\infty}$ acts on $\Lambda$ through a finite quotient.
\endproclaim
\definition{Proof} Because $\Lambda$ is discrete and finitely generated there is an open disc $D$ around $0\in F_{\infty}$ such that the finite set $D\cap\Lambda$ generates $\Lambda$ as an $A$-module. The absolute value on $\overline F_{\infty}$ is invariant with respect to the Galois action therefore the action of Gal$(\overline F_{\infty}|F_{\infty})$ leaves $D\cap\Lambda$ invariant. Every element of Gal$(\overline F_{\infty}|F_{\infty})$ fixing $D\cap\Lambda$ must act trivially on $\Lambda$, too, therefore the action of Gal$(\overline F_{\infty}|F_{\infty})$ on $\Lambda$ has finite image.\ $\square$
\enddefinition
Let $B$ denote the group scheme of invertible upper triangular two by two matrices.
\proclaim{Lemma 8.8} Every finite subgroup of $GL_2(A)$ is either conjugate to a subgroup of $GL_2(\Bbb F_2)\subset GL_2(A)$ or it is a conjugate to a subgroup of $B(A)\subset GL_2(A)$.
\endproclaim
\definition{Proof} Let $G$ be a finite subgroup of $GL_2(A)$. Because $GL_2(A)$ acts without inversions on the Bruhat-Tits tree $\Cal T$ via its natural embedding into $GL_2(F_{\infty})$, so does $G$. Every finite group acting on a tree without inversion has a fixed vertex. Therefore $G$ has a fixed vertex with respect to the action above. Hence it will be conjugate to a subgroup of the stabilizer of the vertex $v_n$ for some $n$. The claim now follows from Proposition 3 of 1.6 in [28] on pages 86-87.\ $\square$
\enddefinition
\proclaim{Proposition 8.9} Let $\phi:A\rightarrow F_{\infty}\{\tau\}$ be a Drinfeld module of rank two over $F_{\infty}$ and fix a
$\bold f_{\goth p}$-linear isomorphism $i:\phi[\goth p]
\rightarrow\bold f_{\goth p}^2$. Then the image of the homomorphism $h:
\text{\rm Gal}(\overline F_{\infty}|F_{\infty})\rightarrow GL_2(\bold f_{\goth p})$ induced by this identification is conjugate to a subgroup of the reduction of a subgroup $G\subset GL_2(A)$ modulo $\goth p$, where $G$ is either $GL_2(\Bbb F_2)$ or $B(A)$.
\endproclaim
\definition{Proof} By Theorem 4.6.9 of [12] on pages 76-78 there is an 
$F$-lattice $\Lambda\subset\overline F_{\infty}$, which has rank two as an $A$-module, such that the Drinfeld module $\phi$ is uniformized by the Carlitz exponential $e_{\Lambda}$ (defined in 4.2.3 of [12] on page 65). In particular for every ideal $\goth n\triangleleft A$ the Galois module $\goth n^{-1}\Lambda/\Lambda$ is isomorphic to $\phi[\goth n]$ where $\goth n^{-1}\subseteq F$ is the fractional ideal: $\{a\in F|ab\in A(\forall b\in\goth n)\}$. The claim is now obvious from Lemmas 8.7 and 8.8.\ $\square$
\enddefinition
 Assume now that $q=2$.
\proclaim{Theorem 8.10} Let $\phi:A\rightarrow F\{\tau\}$ be a Drinfeld module of rank two which has good reduction at every place of $F$ different from $\infty$. Assume that $\deg(\goth p)\geq3$. Then $\phi$ does not have a $\text{\rm Gal}(\overline F|F)$-invariant cyclic $\goth p$-torsion subgroup.
\endproclaim
\definition{Proof} Assume that the claim is false. Then there is a Drinfeld module $\phi$ which satisfies the condition in the claim above such that there is an exact sequence of Gal$(\overline F|F)$-modules:
$$0@>>>V_0@>>>\phi[\goth p]@>>>V_1@>>>0$$
such that the Galois modules $V_0$ and $V_1$ have dimension one as vector spaces over $\bold f_{\goth p}$. Let $\rho_i:\text{Gal}(\overline F|F)\rightarrow\bold f_{\goth p}^*$ denote the modular Galois representation corresponding to the Gal$(\overline F|F)$-module $V_i$ (where $i=1$ or $2$). For every place $v$ of $F$ choose a decomposition group $D_v<\text{Gal}(\overline F|F)$. The restriction of $\rho_i$ to $D_v$ is unramified when $v$ is different from $\goth p$ or $\infty$ since the Drinfeld module $\phi$ has good reduction at every place of $F$ different from $\infty$. We also
know that the restriction of $\rho_i$ to $D_{\infty}$ is at most tamely
ramified. Because the residue field of $F_{\infty}$ has trivial multiplicative group every tamely ramified abelian representation of $D_{\infty}$ is in fact unramified by local class field theory. Hence
$\rho_i$ can be ramified only at $\goth p$.

By Lemma 8.4 the Drinfeld module $\phi$ has good ordinary reduction at $\goth p$. Hence  there is an exact sequence of $D_{\goth p}$-modules:
$$0@>>>U_0@>>>\phi[\goth p]@>>>U_1@>>>0$$
such that the Galois modules $U_0$ and $U_1$ have dimension one as vector spaces over $\bold f_{\goth p}$ by Proposition 8.6. Moreover the $D_{\goth p}$-module $U_1$ is unramified by part $(ii)$ of Proposition 8.6. It is also isomorphic to one of the $D_{\goth p}$-modules $V_0$ and $V_1$. Hence either for $i=0$ or for $i=1$ the Galois representation $\rho_i$ is everywhere unramified. Fix such an $i$. Let $\phi'$ denote the Drinfeld module $\phi$ when $i=0$ and let  $\phi'$ denote the Drinfeld module we get from $\phi$ by dividing out by the Galois-invariant $\goth p$-torsion subgroup $V_0$ when $i=1$. In both cases the $\goth p$-torsion of the Drinfeld module $\phi'$ contains a Gal$(\overline F|F)$-submodule isomorphic to $V_i$ and
$\phi'$ has potentially good reduction at every place of $F$ different from $\infty$.
 
By Proposition 8.9 the image of $D_{\infty}$ under $\rho_i$ is a subgroup of either $B(A)$ or $GL_2(\Bbb F_2)$. It is also a subgroup of $\bold f_{\goth p}^*$ hence it is a group of odd order. Every finite subgroup of $B(A)$ is a $2$-group and the finite group $GL_2(\Bbb F_2)$ only has elements of order 1, 2 or 3. Hence the image of $D_{\infty}$ under $\rho_i$ is either trivial or a group of order $3$. The latter is only possible when $3$ divides the order of $\bold f_{\goth p}^*$. The latter holds if and only if $\deg(\goth p)$ is even. For every $n\in\Bbb N$ let $F_n$ denote the constant field extension of $F$ of degree $n$, that is the unique everywhere unramified extension of $F$ of degree $n$. We get that the Drinfeld module $\phi'$ has a non-zero $F_n$-rational $\goth p$-torsion point where $n=1$ if $\deg(\goth p)$ is odd, and $n=3$, otherwise. In the latter case we have $\deg(\goth p)\geq4$. Now it is clear
the theorem follows from the lemma below.\ $\square$
\enddefinition
Let $T$ denote the place of $F$ corresponding to the prime ideal $(T)\triangleleft A$ as well. 
\proclaim{Lemma 8.11} Let $\phi:A\rightarrow F\{\tau\}$ be a Drinfeld module of rank two over $F$ and let $\goth n\triangleleft A$ be a non-zero ideal. Assume that $\phi$ has an $F_n$-rational non-zero $\goth n$-torsion point for some $n<\deg(\goth n)$ which is not an $\goth m$-torsion point for any proper divisor $\goth m$ of $\goth n$. Then $\phi$ does not have potentially good reduction at the place $T$.
\endproclaim
\definition{Proof} Assume that the claim is false and let $\phi$ be a Drinfeld module which is a counterexample. If we set $\phi(T)=T+g\tau+\Delta\tau^2\in F\{\tau\}$ then the $j$-invariant $j(\phi)=g^3/\Delta$ of $\phi$ lies in $\Cal O_T$ by assumption. For every $k\in F_T$ let $\phi_k:A\rightarrow F_T\{\tau\}$ denote the unique Drinfeld module with the property: $\phi_k(T)=T+kg\tau+k^3\Delta\tau^2$. Clearly $\phi_k$ is isomorphic to $\phi$ over $F_T$. For a suitable choice of $k$ we have $0\leq T(k^3\Delta)\leq 2$. Fix such a $k$ and let $h$, $\Gamma\in F_T$ be such that $\phi_k(T)=T+h\tau+\Gamma\tau^2$. Let $K$ denote the unique
unramified extension of $F_T$ of degree $n$. Let $\Cal R$ be the valuation
ring of $K$.

Assume first that $T(\Gamma)=0$. Because $\phi_k$ is isomorphic to the base change of $\phi$ to $F_T$ there is a $\goth n$-torsion point $y\in(^{\phi_k}K)_{tors}$ which is not an $\goth m$-torsion point for any proper divisor $\goth m$ of $\goth n$. By assumption $T(\Delta)\leq 3T(g)$ hence we have $T(h)=T(kg)\geq0$. Therefore $\phi_k$ gives rise to a Drinfeld module over $\Cal O_T$. In particular $y\in\Cal R$ and the reduction map modulo the maximal ideal of $\Cal R$ is injective restricted to the $A$-module $\{\phi_k(t)(y)|t\in A\}$. The latter has order $2^{\deg(\goth n)}$ but the cardinality of the residue field of $K$ is $2^n$. This is a contradiction.

Assume now that $T(\Gamma)\neq0$. The separable extension $L$ of $K$ we get by adjoining a third root $\lambda$ of the element $\Gamma$ has degree $3$ over $K$. This is a ramified extension because $3$ does not divide the valuation of $\Gamma$ by assumption. Hence the residue field of $L$ is the same as the residue field of $K$. Let $\psi:A\rightarrow L\{\tau\}$ denote the Drinfeld module of rank two such that $\psi(T)=T+\lambda^{-1}h\tau+\tau^2$. Then $\psi$ is isomorphic to the base change of $\phi$ to $L$ hence there is a $\goth n$-torsion point $y\in(^{\psi}L)_{tors}$ which is not an $\goth m$-torsion point for any proper divisor $\goth m$ of $\goth n$. On the other hand $\psi(T)\in\Cal O_L\{\tau\}$, where $\Cal O_L$ is the valuation ring of $L$, hence $\psi$ gives rise to a Drinfeld module over $\Cal O_L$. Now we may conclude the argument as above.\ $\square$
\enddefinition
\proclaim{Proposition 8.12} Let $\phi:A\rightarrow F\{\tau\}$ be a Drinfeld module of rank two which has j-invariant $j(\phi)=0$. Assume that $\deg(\goth p)$ is even and at least $4$. Then $\phi$ does not have a $\text{\rm Gal}(\overline F|F)$-invariant cyclic $\goth p$-torsion subgroup.
\endproclaim
\definition{Proof} By assumption $\phi(T)=T+k\tau^2$ for some $k\in F^*$. Let $\psi:A\rightarrow F\{\tau\}$ denote the unique Drinfeld module with the property: $\psi(T)=T+\tau^2$. Let $K=\Bbb F_4(T)$ be the unique everywhere unramified extension of $F$ of degree two and let $B=\Bbb F_4[T]$ be the integral closure of $A$ in $K$. Both $\phi$ and $\psi$ can be extended uniquely to an $\Bbb F_4$-algebra homomorphism $B\rightarrow K\{\tau\}$ which will be denoted by the same symbol. These homomorphisms are Drinfeld $B$-modules of rank one over $K$. If $l\in\overline F{}^*$ is a third root of $k$ then the linear map $x\mapsto lx$ induces an isomorphism of Drinfeld $B$-modules from $\phi$ to $\psi$ over $\overline F$. Now assume that the claim of the proposition is false and let $M\subset\phi[\goth p](\overline F)$ be a $\text{\rm Gal}(\overline F|F)$-invariant $A$-submodule which has dimension one as a vector space over $\bold f_{\goth p}$. Let $M'=lM$ denote the image of $M$ under the linear map $x\mapsto lx$. For every $\sigma\in\text{Gal}(\overline F|F)$ and $x\in\overline F$ let $x^{\sigma}$ denote the image of $x$ under
$\sigma$. Because $l^{\sigma}/l\in\Bbb F_4^*$ for every $\sigma\in \text{Gal}(\overline F|F)$ we get that
$$\bigcup_{\sigma\in\text{Gal}(\overline F|F)}
\!\!\!\!\!\!\!\{x^{\sigma}|x\in M'\}\subseteq
\bigcup_{\epsilon\in\Bbb F^*_4}\epsilon M'.$$
This property can be reformulated as follows. Let $\Bbb P$ denote the set of $\bold f_{\goth p}$-linear subspaces of $\psi[\goth p]
(\overline F)$ of dimension one. Because the action of $\text{Gal}(\overline F|F)$ on $\psi[\goth p](\overline F)$ is $\bold f_{\goth p}$-linear there is an induced action on $\Bbb P$. By the above the $\text{Gal}(\overline F|F)$-orbit of $M'$ in $\Bbb P$ lies in its orbit with respect to the action of $\Bbb F^*_4$. We will show that $\Bbb P$ has no such element. Because $\deg(\goth p)$ is even the prime ideal $\goth p$ splits to the product of two different prime ideals $\goth p_0\triangleleft B$ and $\goth p_1\triangleleft B$ in $B$. Clearly $\psi[\goth p](\overline F)=\psi[\goth p_0](\overline F)\oplus\psi[\goth p_1](\overline F)$ and for $i=0$, $1$ the $B$-module $\psi[\goth p_i](\overline F)$ has dimension one as a vector space over $\bold f_{\goth p}$. Moreover $\psi[\goth p_i](\overline F)$ is invariant under the action of Gal$(\overline F|K)$. Let $\epsilon_i:\text{Gal}(\overline F|K)\rightarrow\bold f_{\goth p}^*$ denote the modular Galois representation corresponding to the Gal$(\overline F|K)$-module $\psi[\goth p_i](\overline F)$ (where $i=0$ or $1$). The  Galois representation $\epsilon_i$ is unramified at the place $\goth p_{1-i}$ because the Drinfeld module $\psi$ has good reduction at every place of $K$ different from the unique place above $\infty$.

We claim that the restriction of $\epsilon_i$ to the inertia group at $\goth p_i$ is a character of order $2^{\deg(\goth p)}-1$. Let $f_i\in\goth p_i$ be a monic generator. Then
$$\psi(f_i)=\sum_{j=0}^{\deg(\goth p)/2}a_{ij}\tau^j\in B\{\tau\}.$$
We have $a_{i0}=f_i$ and $a_{i\deg(\goth p)/2}=1$ because $f_i$ is monic. Moreover $a_{ij}\in\goth p_i$ for $j=0,1,\ldots,\deg(\goth p)/2-1$ because the base change of the finite, flat group scheme $\psi[\goth p_i]$ to the residue field of $K_i$ is connected where $K_i$ denotes the completion of $K$ with respect to $\goth p_i$. Note that $t\in\overline K_i$ is a non-zero element of $\psi[\goth p_i](\overline K_i)$ if and only if it is a root of the polynomial:
$$g_i(x)=a_{i0}+a_{i1}x^3+\cdots+a_{ik}x^{4^k-1}+\cdots+
a_{i\deg(\goth p)/2}x^{2^{\deg(\goth p)}-1}.$$
Therefore it will be sufficient to prove that $g_i$ is irreducible in
$\widetilde K_i[x]$ where $\widetilde K_i$ is the maximal unramified extension of $K_i$. The latter follows at once from the Sch\"onemann-Eisenstein criterion. 

Therefore the order of the character
$\epsilon_0\cdot\epsilon_1^{-1}$ is $2^{\deg(\goth p)}-1$ because its order is already $2^{\deg(\goth p)}-1$ restricted to the inertia group at $\goth p_i$ (where $i=0$ or $1$). Hence the action of Gal$(\overline F|K)$ on the complement of $\{\psi[\goth p_0](\overline F)\}\cup\{\psi[\goth p_1](\overline F)\}$ in $\Bbb P$ is transitive. Because $2^{\deg(\goth p)}-1>3$ we get that no element of this complement can have the property above. Let $\sigma\in\text{Gal}(\overline F|F)$ be an element whose image in the quotient group $\text{Gal}(K|F)$ is the generator. The action of $\sigma$ on the places of $K$ will interchange $\goth p_0$ and $\goth p_1$ so its action on $\text{Gal}(\overline F|K)$ via conjugation interchanges the characters
$\epsilon_0$ and $\epsilon_1$. Hence it must interchange the submodules $\psi[\goth p_0](\overline F)$ and $\psi[\goth p_1]
(\overline F)$ as well. Since these elements of $\Bbb P$ are fixed by the action of $\Bbb F_4^*$ they could not have the property either.\ $\square$
\enddefinition

\heading 9. The proofs of Theorem 1.2 and Corollary 1.7
\endheading

\definition{Notation 9.1} Let $\goth m\triangleleft A$ be an arbitrary non-zero ideal. For every $(\alpha,\beta)\in(A/\goth m)^2$, and $n$ integer let
$$W_{\goth m}(\alpha,\beta,n)=\{0\neq (a,b)\in A^2|
(a,b)\equiv(\alpha,\beta)\!\!\!\!\mod\goth m,-n=
\min(\infty(a),\infty(b))\}.$$
For every $(\alpha,\beta)$ as above and $N$ positive integer let $\epsilon_{\goth m}(\alpha,\beta,N)(z)$ denote the function:
$$\epsilon_{\goth m}(\alpha,\beta,N)(z)=
\prod_{n\leq N}\left(
\prod_{(a,b)\in W_{\goth m}(\alpha,\beta,n)}
\!\!\!\!\!\!\!(az+b)\cdot\!\!\!\!\!\!\!\!\!\!\!
\prod_{(c,d)\in W_{\goth m}(0,0,n)}\!\!\!\!\!\!\!
(cz+d)^{-1}\right).$$
on the set $\Omega$. The latter is clearly holomorphic in the variable $z$.  According to Lemma 4.5 of [17] on pages 145-146 the limit
$$\epsilon_{\goth m}(\alpha,\beta)(z)=\lim_{N@>>>\infty}
\epsilon_{\goth m}(\alpha,\beta,N)(z)$$
converges uniformly in $z$ on every admissible open subdomain of $\Omega$ and defines a ho\-lo\-mor\-phic function.
\enddefinition
\definition{Definition 9.2} If $\psi:A@>>>\Bbb
C_{\infty}\{\tau\}$ is a Drinfeld module of rank two over $A$, then 
$$\psi(T)=T+g(\psi)\tau+\Delta(\psi)\tau^2,$$
where $g$ and $\Delta$ are Drinfeld modular forms of weight $q-1$ and $q^2-1$, respectively. Under the identification of Theorem 7.2 $g$ and $\Delta$ are holomorphic functions on $\Omega$ and they are equal to:
$$g(z)=\!\!\!\!\!\!\!\!\!\!\!\!\!\!\!\!\!\!\!\!
\sum_{(0,0)\neq(\alpha,\beta)\in(A/(T))^2}
\!\!\!\!\!\!\!\!\!\!\!\!\!\!\!\!\!\!\!
\epsilon_{(T)}(\alpha,\beta)(z),\quad\Delta(z)=\!\!\!\!\!\!\!\!\!\!\!\!\!\!\!\!\!\!\!\!
\prod_{(0,0)\neq(\alpha,\beta)\in(A/(T))^2}
\!\!\!\!\!\!\!\!\!\!\!\!\!\!\!\!\!\!\!
\epsilon_{(T)}(\alpha,\beta)(z),$$
which is an immediate consequence of the uniformization theory of
Drinfeld modules over $\Bbb C_{\infty}$.
\enddefinition
Let $K$ be an extension of $F$ and let $P$ be a $K$-valued point of $X_0(\goth p)$. Let $(\phi,G)$ be the Drinfeld module of rank two over $K$ equipped with a Gal$(\overline K|K)$-invariant cyclic $\goth p$-torsion subgroup $G$ which corresponds to the point $P$. The $j$-invariant $j(P)$ of $P$ is the $j$-invariant $j(\phi)\in K$ of the Drinfeld module $\phi$ by definition.
\proclaim{Lemma 9.3} Assume that $q=2$. If the reduction of the point $P\in Y_0(\goth p)(F_{\infty})$ lies on the irreducible component
$\widetilde w$ then $j(P)\in\Cal O_{\infty}$.
\endproclaim
\definition{Proof} Let $Q$ be the image of $P$ with respect to the map 
$Y_0(\goth p)\rightarrow Y_0(1)$ corresponding to forgetting the level structure. By definition $j(P)$ is just the image of $Q$ under the isomorphism $Y_0(1)\cong\Bbb A^1$ induced by the $j$-invariant. According to Theorem 7.2 the rigid-analytic space underlying
$Y_0(1)$ is isomorphic to the generic fiber of the formal scheme $GL_2(A)\backslash\widehat\Omega$. The incidence graph of the special fiber of $GL_2(A)\backslash\widehat\Omega$ is canonically isomorphic to the quotient graph $GL_2(A)\backslash\Cal T$. Because the vertex $w$ of $\Cal T$ is equivalent to $v_0$ modulo the action of $GL_2(A)$ the reduction of $Q$ lies on the irreducible component of the special fiber of $GL_2(A)\backslash\widehat\Omega$ which is the image of $v_0$ under the quotient map $\Cal T\rightarrow GL_2(A)\backslash\Cal T$. Hence by Corollary 1.9 of [6] on page 186 there is a $z\in V$ whose image with respect to the uniformization in Theorem 7.2 is the point $Q$ where $V$ denotes the rational subdomain whose underlying set of points is:
$$V=\{z\in\Bbb C_{\infty}|\ |z-c|=1\ (\forall c\in\bold f_{\infty})\}\subset\Omega.$$
Hence it will be sufficient to show that $|j(z)|\leq1$ for every $z\in V$ where $j$ is the Drinfeld modular function $j:\Omega\rightarrow\Bbb C_{\infty}$. Let $z_0\in V$ be an element of $\Bbb F_4-\Bbb F_2$. Then $A+Az_0=\Bbb F_4[T]$ is a rank two discrete $A$-lattice which corresponds to a Drinfeld module $\psi:A\rightarrow\Bbb C_{\infty}\{\tau\}$ of rank two by the uniformization theory of Drinfeld modules. As explained in the proof of Theorem 2.13 of [6] on page 192 the $j$-invariant of $\psi$ is zero. Hence $\psi(T)=\overline\pi\tau^2+T$ for some $\overline\pi\in\Bbb C^*_{\infty}$. The modular unit $\epsilon_{(T)}(\alpha,\beta)(z_0)$ is a root of the additive polynomial $\overline\pi x^4+Tx$ corresponding to $\psi(T)$ for every $(\alpha,\beta)\in(A/(T))^2$. Therefore $|\epsilon_{(T)}(\alpha,\beta)(z_0)|=|T/\overline\pi|^{1/3}$ for every non-zero $(\alpha,\beta)$ as above. The reduction of the affinoid $V$ is irreducible hence every nowhere zero holomorphic function on $V$ has constant absolute value. In particular $|\epsilon_{(T)}(\alpha,\beta)(z)|=|T/\overline\pi|^{1/3}$ for every $z\in V$ and for every non-zero $(\alpha,\beta)$ as above. We get that $|\Delta(z)|=|T/\overline\pi|$ and $|g(z)^3|\leq|T/\overline\pi|$ for every $z\in V$ using the ultrametric inequality. The claim is now clear.\ $\square$
\enddefinition
\proclaim{Proposition 9.4} For every place $v$ of $F$ the reduction map $W(\goth p)(F)\rightarrow\Cal W(\goth p)(\bold f_v)$ is injective.
\endproclaim
\definition{Proof} By Lemma 3.4 the order of the finite group $W(\goth p)(F)$ is not divisible by $p$. Hence the claim is obvious for every place $v$ where $W(\goth p)$ has good reduction. Because $W(\goth p)$ is a quotient of $J_0(\goth p)$ this holds for every place $v$ different from $\goth p$ and $\infty$. At the places $\goth p$ and
$\infty$ the Jacobian $J_0(\goth p)$ has multiplicative reduction, that is, the connected component of the identity in the fiber of
its N\'eron model over these closed points is a torus. Therefore $W(\goth p)$ also has multiplicative reduction at these places. The claim now follows from part $(i)$ of Lemma 7.13 of [17] on page 162 in these cases.\ $\square$
\enddefinition
\definition{Notation 9.5} By Proposition 9.3 of [2] on pages 586-587 we know that for every non-zero ideal $\goth n\triangleleft A$ the projective curve $X_0(\goth n)$ has a smooth, projective model $\overline M_0(\goth n)$ over $U_{\goth n}$ where $U_{\goth n}$ is the
complement of the support of the ideal $\goth n$ in Spec$(A)$. In
the paper quoted above it is also proved that the Zariski closure of
the cusps in $\overline M_0(\goth n)$ is a finite, \'etale scheme
over $U_{\goth n}$. When $\goth n=\goth p$ is a prime ideal Gekeler constructed a model of $X_0(\goth p)$ over Spec$(A)$ which contains $\overline M_0(\goth p)$ as an open subscheme. We are going to denote the latter by $\overline M_0(\goth p)$ as well by slight abuse of notation. For the relevant properties of this scheme see 5.1-5.8 of [5], pages 229-233.
\enddefinition
\definition{Proof of Theorem 1.2} We are going to use again the notation which we used in chapter 7. Let $w_{\goth p}$ denote again the matrix which was introduced in the proof of Lemma 6.6. As we already mentioned the latter normalizes the group $\Gamma$ hence its M\"obius action on $\Omega$ descends to the quotient $\Gamma\backslash\Omega$. Under the identification of Theorem 7.2 the corresponding automorphism of $Y_0(\goth p)$ is the Atkin-Lehner involution. Because $\Gamma_{tor}$ is a characteristic subgroup of $\Gamma$ the matrix $w_{\goth p}$ normalizes this subgroup of $GL_2(F_{\infty})$, too. Hence its M\"obius action on $\Omega$ descends to the quotient $\Gamma_{tor}\backslash\Omega$ as well. This action extends to $\Sigma$ uniquely and leaves the cover $\Cal V$ invariant. 
The corresponding automorphism of the formal scheme $\widehat\Sigma$ commutes with the action of $\Gamma'$ hence we get an involutive automorphism of the Spec$(\Cal O_{\infty})$-scheme $\goth X_0(\goth p)$ by formal GAGA. We will call the latter the Atkin-Lehner 
involution of the scheme $\goth X_0(\goth p)$ and will denote it by
$W_{\goth p}$ by the usual abuse of notation. This is justified as its action of the generic fiber of $\goth X_0(\goth p)$ is the usual Atkin-Lehner involution by the above. Moreover the action on the incidence graph $\Cal G(\goth p)$ of the special fiber of $\goth X_0(\goth p)$ induced by $W_{\goth p}$ is the Atkin-Lehner involution which  we introduced in the proof of Lemma 6.6.

Let $P$ be an $F$-rational point of $X_0(\goth p)$. By Proposition 7.12 the model $\goth X_0(\goth p)$ is regular hence the reduction of $P$ lies in $\goth X^*_0(\goth p)(\bold f_{\infty})$. The ordinary double points of the special fiber of $\goth X_0(\goth p)$ are $\bold f_{\infty}$-rational. Hence if $v$ is an irreducible component of the special fiber of $\goth X_0(\goth p)$ whose degree, considered as a vertex of the incidence graph, is three, then every $\bold f_{\infty}$-rational point of $v$ is an ordinary double point. Therefore the reduction of $P$ must lie on the irreducible component $\widetilde u_{-1}$, or its image with respect to the Atkin-Lehner involution, or the irreducible component
$\widetilde w$ by Propositions 6.11 and 6.13. The first two irreducible components both have degree two as a vertex of $\Cal G(\goth p)$, hence they have only one $\bold f_{\infty}$-rational point which is not an ordinary double point. These must be the reductions of the cusps $0$ and $\infty$, respectively. Therefore the reduction of $P$ is either equal to the reduction of one of the cusps or $\deg(\goth p)$ is even and the reduction of $P$ lies on the irreducible component $\widetilde w$.

Assume first that the reduction of $P$ in $\goth X^*_0(\goth p)(\bold f_{\infty})$ is the same as the reduction of some cusp. By applying the Atkin-Lehner involution to $P$, if necessary, we may assume that this cusp is $0$. In this case the image of the divisor class $P-0$ in $W(\goth p)(F)$ with respect to the quotient map is equal to the image of the trivial divisor class in $W(\goth p)(F)$ by Proposition 9.4. Therefore we must have $P=0$ by the formal immersion property (Proposition 7.14). Hence we may assume that $\deg(\goth p)$ is even and the reduction of $P$ lies on the irreducible component $\widetilde w$. By Lemma 9.3 above we know that $j(P)\in\Cal O_{\infty}$.

We are going to show that $j(P)\in\Cal O_v$ for every other place
$v$ of $F$. Assume that the claim is false. Because the model
$\overline M_0(\goth p)$ is projective over Spec$(A)$ every $F$-valued point on $X_0(\goth p)$ has a unique extension to a section Spec$(A)\rightarrow\overline M_0(\goth p)$. By assumption the section extending $P$ must intersect the section extending one of the cusps in the fiber of $\overline M_0(\goth p)$ at $v$. We may assume that $P$ intersects the section extending $\infty$ as the other case can be treated exactly the same way. Then the image of the linear equivalence class of the divisor $P-\infty$ with respect to the reduction map $J_0(\goth p)(F)\rightarrow\Cal J_0(\goth p)(\bold f_v)$ is zero. For every $F$-rational divisor $D$ of degree zero on $X_0(\goth p)$ let $[D]$ denote the image of the linear equivalence class of $D$ with respect to the quotient map $J_0(\goth p)\rightarrow W(\goth p)$. By the universal property of the N\'eron model we get that 
the image of $[P-\infty]$ with respect to the reduction map $W(\goth p)(F)\rightarrow\Cal W(\goth p)(\bold f_v)$ is zero. Hence $[P-\infty]$ is zero in $W(\goth p)(F)$ by Proposition 9.4. By claim $(iii)$ of Proposition 6.11 the reduction of $W_{\goth p}(P)$ lies on the same irreducible component as the reduction of $P$. Since this irreducible component is a rational curve which has only one ordinary double point the difference of these reductions is linearly equivalent to zero on the special fiber of $\goth X_0(\goth p)$. Hence the image of the linear equivalence class of the divisor $P-W_{\goth p}(P)$ with respect to the reduction map $J_0(\goth p)(F)\rightarrow\Cal J_0(\goth p)(\bold f_{\infty})$ is zero. Using the same argument as above we get from this fact that $[P-W_{\goth p}(P)]$ is zero in $W(\goth p)(F)$ by Proposition 9.4. Therefore we may compute as follows:
$$\split 0=W_{\goth p}[P-\infty]=[W_{\goth p}(P)-W_{\goth p}(\infty)]=&
[W_{\goth p}(P)-P]+[P-W_{\goth p}(\infty)]\\
=&[P-0]=[P-\infty]+[\infty-0]=[\infty-0]\endsplit$$
where $W_{\goth p}$ denotes the endomorphism of $W(\goth p)$ induced by the Atkin-Lehner involution, too. The linear equivalence class of the divisor $\infty-0$ generates the group $\Cal C(\goth p)$ hence the image of the latter with respect to the map $J_0(\goth p)\rightarrow W(\goth p)$ is zero. This is a contradiction by claim $(ii)$ of Proposition 5.12. We conclude that $j(P)\in\Cal O_v$ for every place $v$ of $F$. This is only possible if $j(P)\in\Bbb F_2$. By Proposition 8.12 we may assume that $j(P)=1$. Every Drinfeld module $\phi$ of rank two over $F$ with $j$-invariant $j(\phi)=1$ is isomorphic to the Drinfeld module $\psi:A\rightarrow F
\{\tau\}$ such that $\psi(T)=T+\tau+\tau^2$ over $F$ because the latter has no non-trivial automorphisms over $\overline F$. The Drinfeld module $\psi$ has good reduction at every place $v\neq\infty$ of $F$. This is not possible by Theorem 8.10.\ $\square$
\enddefinition
\definition{Proof of Corollary 1.7} The scheme $X_0(\goth p)$ is smooth over Spec$(F)$. Hence the scheme
$\bold{Isom}_{\text{Spec}(F)}(X_0(\goth p),X_0(\goth p))$ representing the functor associating to each Spec$(F)$-scheme $S$ the set of $S$-automorphisms of the scheme $X_0(\goth p)\times_{\text{Spec}(F)}S$ is finite and unramified over Spec$(F)$ by Theorem 1.11 of [1] on pages 84-85. Therefore every automorphism $u$ of $X_0(\goth p)$ is actually defined over the separable closure of $F$. Let
Aut$(J_0(\goth p)_{\overline F})$ denote the group of automorphism of the group scheme $J_0(\goth p)$ over $\overline F$. The homomorphism from Aut$(X_0(\goth p))$ into Aut$(J_0(\goth p)_{\overline F})$ induced by Picard functoriality is injective. Moreover this map is also Gal$(\overline F|F)$-invariant. Hence every automorphism $u$ of $X_0(\goth p)$ is actually defined over $F$ by claim $(i)$ of Lemma 5.6. Therefore the image of the cusps of $X_0(\goth p)$ under $u$ are $F$-rational points. By Theorem 1.2 the set $X_0(\goth p)(F)$ consists of the cusps hence every automorphism of $X_0(\goth p)$ fixes the two cusps or interchanges them. The Atkin-Lehner involution is an automorphism which interchanges the cusps. Therefore it will be sufficient to prove that every involution $u$ of the curve $X_0(\goth p)$ fixing the cusps is the identity by claim $(iii)$ of Lemma 5.6. Let $f$ be a non-zero rational function on $X_0(\goth p)$ whose divisor is $N(\goth p)(\infty-0)$. Because the automorphism $u$ fixes the points $0$ and $\infty$ the pull-back $u^*(f)$ of $f$ with respect to $u$ is also a non-zero rational function on $X_0(\goth p)$ whose divisor is $N(\goth p)(\infty-0)$. Such a function is unique up to a non-zero scalar, so $u^*(f)=\alpha f$ for some $\alpha\in F^*$. But $u$ is an involution so $f=u^*(u^*(f))=\alpha^2f$ hence $\alpha^2=1$. The field $F$ has characteristic $2$ so the equation $\alpha=1$ must hold. Hence the rational function $f$ is left invariant by the involution $u$. Assume now that $u$ is not the identity. Let $X_0(\goth p)^u$ be the quotient of $X_0(\goth p)$ by $u$ and let $\pi:X_0(\goth p)\rightarrow X_0(\goth p)^u$ denote the quotient map. By the above there is map $g:X_0(\goth p)^u\rightarrow\Bbb P^1$ such that $f=g\circ\pi$. The degree of $f$ is $N(\goth p)$ which is an odd number. On the other hand the degree of $\pi$ is two. This is a contradiction.\ $\square$
\enddefinition
Actually our argument works for every field $\Bbb F_q(T )$ of characteristic $2$ under the condition that $Y_0(\goth p)$ has no $\Bbb F_q(T )$-rational points.
\definition{Remark 9.6} Let $B$ denote again the group scheme of invertible upper triangular two by two matrices. Let $h:\Gamma_0(\goth p)\rightarrow\Bbb Z/N(\goth p)\Bbb Z$ denote the composition of the reduction map $r:\Gamma_0(\goth p)\rightarrow
 B(\bold A/\goth p)\subset GL_2(\bold A/\goth p)$ mod $\goth p$, the upper left conner element $a:B(\bold A/\goth p)\rightarrow(\bold A/\goth p)^*$ and the surjection $s:(\bold A/\goth p)^*\rightarrow\Bbb Z/N(\goth p)\Bbb Z$ which is unique up to isomorphism. We used that this homomorphism factors through $\overline{\Gamma}_0(\goth p)$ at a crucial point both in [17] and [18], but did not give a proof of this not entirely obvious fact. This small gap is easy to fill with some of the methods of this paper. Because $\overline{\Gamma}_0(\goth p)=(\Gamma_0(\goth p)_{tors})^{ab}$ it will be sufficient to show that $h$ is trivial restricted to every finite subgroup $G$ of $\Gamma_0(\goth p)$. As we saw in the proof of Lemma 8.8 the action of $G$ on the Bruhat-Tits tree has a fixed vertex. Hence by the results of [7] it will be conjugate to a subgroup of $B(A)$ or to a subgroup of $GL_2(\Bbb F_q)$, where the second case occurs only when $\deg(\goth p)$ is even. In any case every element $g$ of $G$ has order dividing $q^k(q-1)$ for some $k\in\Bbb N$, if $\deg(\goth p)$ is odd, and has order dividing $q^k(q^2-1)$ for some $k\in\Bbb N$, otherwise. Hence $a(r(g))$ is mapped to the identity by $s$. 
\enddefinition 

\heading 10. The proof of Theorem 1.4
\endheading

In the next claim and its proof we are not going to assume that $q=2$.
\proclaim{Proposition 10.1} Assume that the following holds:
\roster
\item "$(i)$" The curve $Y_1(\goth n)$ has no $F$-rational points if $\deg(\goth n)\geq3$,
\item "$(ii)$" The curves $Y_0(T^3)$ and $Y_0(T^2(T+1))$ have no rational points.
\endroster
Then the claim of Theorem 1.5 holds for $q$.
\endproclaim
The proposition above is just a slight variation of Schweizer's
Lemma 4.5 of [26] on page 613.
\definition{Proof} Let $\phi:A\rightarrow F\{\tau\}$ be a
Drinfeld module of rank two over $F$ and write
$(^{\phi}F)_{tors}=A/\goth m\oplus A/\goth n$ where $\goth m$
divides $\goth n$. By assumption the latter can be written as $\goth
n=\goth m\goth p$ for some non-zero ideal $\goth p\triangleleft A$.
Assume that $\phi$ is a counterexample to the claim of Theorem 1.2.
By claim $(i)$ above the degree of the ideal $\goth n$ is at most
two. Hence $\deg(\goth m)$ is at least 1. On the other hand we have
$\deg(\goth m)\leq1$ by Proposition 4.4 of [26] on page 613. Let
$v:\phi\rightarrow\psi$ be the cyclic $\goth m$-isogeny whose kernel
is the first component of the direct sum above. According to Lemma
4.1 of [26] on page 611, the Drinfeld module $\psi$ over $F$ has a
cyclic $\goth{mn}$-isogeny defined over $F$. This Drinfeld module
and isogeny give rise to an $F$-rational point on the Drinfeld
modular curve $Y_0(\goth p\goth m^2)$. If the ideals $\goth p$ and
$\goth m$ are different, then there is an automorphism of $A$ which
takes $\goth p$ and $\goth m$ to $(T+1)$ and $(T)$, respectively. On
the other hand if ideals $\goth p$ and $\goth m$ are equal, then
there is an automorphism of $A$ which takes $\goth p=\goth m$ to
$(T)$. The pull-back of $Y_0(\goth p\goth m^2)$ with respect to
these automorphisms is either $Y_0(T^2(T+1))$ or $Y_0(T^3)$,
respectively. But these curves have no rational points by claim
$(ii)$ above which is a contradiction.\ $\square$
\enddefinition
\definition{Remark 10.2} In the course of the proof of Theorem 1.5 we
may restrict to those $\goth n$ ideals which has no proper divisors
of degree at least three while we verify the validity of condition
$(i)$ of Proposition 10.1. In particular we may assume that either
$\goth n$ is a prime ideal of degree at least three or it is the
product of irreducible polynomials of degree at most two. When $q=2$
the only irreducible polynomials of degree $1$ and $2$ in $\Bbb
F_q[T]$ are $T$, $T+1$ and $T^2+T+1$. Moreover there is an
automorphism of $\Bbb F_2[T]$ which maps $T$ to $T+1$. Therefore by
Theorem 1.2 and the claim above in order to prove Theorem 1.5 we
have to show that the following curves have no rational points:
$$Y_0(T(T^2+T+1)),Y_0(T^3),Y_0(T^2(T+1)),Y_1(T^4+T^2+1)$$
taking into account that if the curve $Y_0(\goth n)$ has no rational
points then so does $Y_1(\goth n)$. The first three cases were
settled by Schweizer (Lemma 1.3 of [26] on pages 605-606).
\enddefinition
\definition{Notation 10.3} In the rest of this section we will occupy ourselves
with proving that the curve $Y_1(T^4+T^2+1)$ has no rational points. First we
will collect some facts about the geometry of the curve $X_0(T^4+T^2+1)$. Let
$W_{\goth n}$ denote the full Atkin-Lehner involution of the modular curve
$Y_0(\goth n)$ for every non-zero ideal $\goth n\triangleleft A$. This map
induces an involution of the curve $X_0(\goth n)$ which will be denoted by the
same symbol. The index $\goth n$ will be dropped from the notation
$W_{\goth n}$ when the ideal $\goth n$ can be clearly identified from the
context. Let $X_+(\goth n)$ denote the quotient of the curve $X_0(\goth n)$ by
the involution $W$. Recall that the set of geometric points of $X_0(\goth n)$
in the complement of $Y_0(\goth n)$ are called its cusps.
\enddefinition
\proclaim{Lemma 10.4} The following holds:
\roster
\item"$(i)$" the curve $X_0(T^4+T^2+1)$ has genus $2$,
\item"$(ii)$" the quotient map $X_0(T^4+T^2+1)\rightarrow X_+(T^4+T^2+1)$
is a hyperelliptic cover which has three ramification points,
\item"$(iii)$" the latter are cusps conjugate over $F$.
\endroster
\endproclaim
\definition{Proof} Claim $(i)$ is a consequence of the formula on page 331 for the genus of Drinfeld modular curves in [25]. By formula {\it b)} of
Proposition 7 on page 336 of the same paper the number of fixed
points of the Atkin-Lehner involution is at most $3$. Now claims
$(ii)$ and $(iii)$ follow from the content of Example 10 of [25] on
page 337 and from part $(b)$ of Proposition 1 of the same reference
on page 333.\ $\square$
\enddefinition
\proclaim{Lemma 10.5} The fiber of the model $\overline
M_0(T^4+T^2+1)$ over every prime ideal $\goth q\neq(T^2+T+1)$ of $A$
is an ordinary curve.
\endproclaim
\definition{Proof} Let $X_{\goth q}$ denote simply the fiber of
$\overline M_0(T^4+T^2+1)$ at the prime ideal $\goth
q\neq(T^2+T+1)$. We know that $X_{\goth q}$ is a smooth, projective
curve over $\bold f_{\goth q}$ whose genus is $2$ by flatness.
Because $\overline M_0(T^4+T^2+1)$ is a relatively minimal model for
$X_0(T^4+T^2+1)$ over $U_{T^2+T+1}$ the Atkin-Lehner involution $W$ extends to an involution on $\overline M_0(T^4+T^2+1)$ whose
restriction to $X_{\goth q}$ is a hyperelliptic involution. The
former fixes the closure $Z$ of the zero-dimensional subvariety of
$X_0(T^4+T^2+1)$ supported on the three cusps left invariant by the
Atkin-Lehner involution in $\overline M_0(T^4+T^2+1)$. As we already noted in Notation 9.5 the scheme $Z$ is finite, \'etale over $U_{T^2+T+1}$  hence the canonical hyperelliptic involution of $X_{\goth q}$ has at least three
different fixed points. Hence it will be sufficient to show that
every smooth, projective curve $X$ of genus $2$ defined over a
characteristic two field $\bold k$ such that the canonical
hyperelliptic cover $\pi:X\rightarrow\Bbb P^1_{\bold k}$ has three
ramification points is ordinary. We may assume without the loss of
generality that $\bold k$ is algebraically closed. Let $A$, $B$ and
$C$ denote the ramification points. Then the linear equivalence
classes of the divisors $2B-2A$ and $2C-2A$ are trivial as they are
the pull-backs of the divisors $\pi(B)-\pi(A)$ and $\pi(C)-\pi(A)$
with respect to $\pi$, respectively. Hence the linear equivalence
classes of $B-A$ and $C-A$ are $2$-torsion. If they generated the
same group then they would be equal, so $B$ would be equal to $C$.
Hence the group generated by the linear equivalence classes of $B-A$
and $C-A$ has order four. Therefore the $2$-torsion of the Jacobian
of $X$ has order four so it is an ordinary abelian variety.\
$\square$
\enddefinition
\proclaim{Proposition 10.6} The following holds:
\roster
\item"$(i)$" the Mordell-Weil group $J_0(T^4+T^2+1)(F)$ of the
Jacobian of the Drinfeld modular curve $X_0(T^4+T^2+1)$ is finite,
\item"$(ii)$" for every prime ideal $\goth q\neq(T^2+T+1)$ of $A$ the reduction map
modulo $\goth q$ on $X_0(T^4+T^2+1)(F)$ is injective.
\endroster
\endproclaim
\definition{Proof}  By the main theorem of [24] on page 509 it is enough to show that the Grothendieck $L$-function of the abelian variety $J_0(T^4+T^2+1)$ does not vanish at $q^{-1}$ in order to prove the first claim. We are going to use the notation which we introduced at the beginning of the proof of Proposition 3.8. Fix a prime $l$ different from $2$ and let $\rho$ be an irreducible component of  the $l$-adic representation $H^1(J_0(T^4+T^2+1)_{\overline F},\overline{\Bbb Q}_l)$ of the absolute Galois group of $F$. The curves $X_0(1)$ and $X_0(T^2+T+1)$ have genus zero hence we get that
the degree of the conductor of $\rho$ is $5$ by repeating the argument in the proof which have just mentioned above. 
Therefore the degree of the Grothendieck $L$-function $L(\rho,t)$ as
a polynomial in $t$ is one. In particular the multiplicity of the
zero at $q^{-1}$ is at most one. The Atkin-Lehner operator acts on
the Jacobian of $X_0(\goth n)$ by Picard functoriality which induces
an action on the cohomology group $H^1(J_0(\goth n)_{\overline
F},\overline{\Bbb Q}_l)$. This action commutes with the action of
the absolute Galois group of $F$ and its restriction to every
irreducible component $\rho$ is multiplication by $\pm1$. Moreover
by the analogue of the classical Atkin-Lehner theory for Drinfeld
modular curves (see Remark 10.7 below) the sign of the functional
equation of $L(\rho,t)$ is the opposite of the sign of this
multiplication.  The $+1$-eigenspace of the action of the
Atkin-Lehner operator on $H^1(J_0(T^4+T^2+1)_{\overline
F},\overline{\Bbb Q}_l)$ is isomorphic to
$H^1(X_+(T^4+T^2+1)_{\overline F},\overline{\Bbb Q}_l)$ which is
zero by claim $(ii)$ of Lemma 10.4. Hence the sign of the functional
equation of $L(\rho,t)$ is $+1$ so the order of vanishing at
$q^{-1}$ is even. Claim $(i)$ is now clear.

Because the curve $X_0(\goth n)$ has a smooth, projective model over
$U_{T^2+T+1}$ the N\'eron model of the Jacobian $J_0(T^4+T^2+1)$ is
an abelian scheme $\Cal J$ over $U_{T^2+T+1}$. It will be sufficient to
prove that for every prime ideal $\goth q\neq(T^2+T+1)$ of $A$ the
reduction map modulo $\goth q$ on the Mordell-Weil group
$J_0(T^4+T^2+1)(F)$ is injective in order to show claim $(ii)$. The
reduction map is injective restricted to the prime-to-two torsion of
$J_0(T^4+T^2+1)(F)$. Hence by claim $(i)$ we only need to show that
it is injective restricted to the $2$-torsion as well. Let
$\widetilde F_{\goth q}$ denote the maximal unramified extension of
$F_{\goth q}$ and let $\widetilde{\Cal O}_{\goth q}$, $\bold
k_{\goth q}$ denote the valuation ring and the residue field of the
local field $\widetilde F_{\goth q}$, respectively. The base change
$\Cal J_{\goth q}$ of $\Cal J$ to the geometric point Spec$(\bold
k_{\goth q})\rightarrow U_{T^2+T+1}$ is ordinary by Lemma 10.5, hence its formal group is multiplicative by Lazard's theorem. The latter
has no non-trivial deformations, so the formal group of the base
change of $\Cal J$ to the spectrum of the local ring
$\widetilde{\Cal O}_{\goth q}$ is also multiplicative. The group of
$\widetilde F_{\goth q}$-valued points of $\Cal J$ which map to the
identity of $\Cal J_{\goth q}$ with respect to the reduction map
into the special fiber is isomorphic to the multiplicative group of
those elements of $\widetilde{\Cal O}_{\goth q}$ which are congruent
to $1$ modulo the maximal ideal of $\widetilde{\Cal O}_{\goth q}$.
Hence the former has no $2$-torsion, so the reduction map is even
injective restricted to $J_0(T^4+T^2+1)(\widetilde F_{\goth
q})[2]$.\ $\square$
\enddefinition
\definition{Remark 10.7} In this paragraph we briefly indicate why the analogue of the classical Atkin-Lehner theory is true for Drinfeld modular curves (or alternatively the reader may consult the proof of Proposition 3.14 of [29] on pages 225-226). By 8.3.8 of [8] on page 81 there is an Eichler-Shimura isomorphism:
$$i_{ES}:H^1(J_0(\goth n)_{\overline F},\overline{\Bbb Q}_l)\rightarrow
 H_!(\Cal T,\overline{\Bbb Q}_l)^{\Gamma_0(\goth n)}\otimes_{\overline{\Bbb Q}_l}
\overline{\Bbb Q}_l^2,$$ where $H_!(\Cal T,\overline{\Bbb
Q}_l)^{\Gamma_0(\goth n)}$ is the space of $\overline{\Bbb
Q}_l$-valued cuspidal cochains on the Bruhat-Tits tree $\Cal T$
invariant with respect to the Hecke congruence group $\Gamma_0(\goth n)$. The Eichler-Shimura isomorphism $i_{ES}$ is equivariant to the usual action of the Hecke algebra on these two vector spaces over
$\overline{\Bbb Q}_l$. Moreover the map $i_{ES}$ is equivariant with
respect to the action of the Atkin-Lehner involution in the sense
that under this isomorphism the operator induced by $W$ on the
cohomology group $H^1(J_0(\goth n)_{\overline F},\overline{\Bbb
Q}_l)$ corresponds to the tensor product of the identity of
$\overline{\Bbb Q}_l^2$ and the operator on $H_!(\Cal
T,\overline{\Bbb Q}_l)^{\Gamma_0(\goth n)}$ corresponding to the
product of the full Atkin-Lehner operator and the local Atkin-Lehner
operator at $\infty$ for automorphic forms of conductor $\goth
n\infty$ where we write the product of divisors multiplicatively.
Since the latter acts as multiplication by $-1$ on the elements of
$H_!(\Cal T,\overline{\Bbb Q}_l)^{\Gamma_0(\goth n)}$ the claim is
clear.
\enddefinition
\proclaim{Theorem 10.8} The Drinfeld modular curve
$Y_1(T^4+T^2+1)$ has no $F$-rational points.
\endproclaim
\definition{Proof} Assume that the claim is false, and let $P\in Y_1(T^4+T^2+1)(F)$ be an $F$-rational point. The point $P$ maps to an $F$-rational point $Q$ on $Y_0(T^4+T^2+1)$ with respect to the map $Y_1(T^4+T^2+1) \rightarrow Y_0(T^4+T^2+1)$ of these coarse moduli induced by the forgetful functor. By claim $(ii)$ of Proposition 10.6 the reduction of $Q$ modulo $T$ is different from the cusps. On the other hand the latter is not possible by Lemma 8.11.\ $\square$
\enddefinition

\enddocument